\documentclass{llncs}
\usepackage{setspace}
\usepackage{amsmath}
\usepackage[boxed]{algorithm2e}
\usepackage{color}
\usepackage{amssymb}
\usepackage[isolatin]{inputenc}
\usepackage[english]{babel}
\usepackage{graphicx}
\usepackage{psfrag}
\usepackage{mathrsfs}
\usepackage[font=footnotesize,format=plain,labelfont=bf,up,textfont=it,up,margin=0.03\textwidth]{caption}
\usepackage[captionskip=12pt]{subfig}
\usepackage{fancyhdr}

\begin{document}
 
\title{A Review of the Finite Cell Method\\for Nonlinear Structural Analysis of Complex CAD and Image-based Geometric Models}
  
  \author{Dominik Schillinger \and Quanji Cai \and Ralf-Peter Mundani \and Ernst Rank}
  
  \institute{Lehrstuhl f\"{u}r Computation in Engineering, Dept. of Civil Engineering and Surveying,\\ Technische Universit\"{a}t M\"{u}nchen, Arcisstr. 21, 80333 M\"{u}nchen, Germany  \\[0.16cm] \email{$\{$schillinger,cai,mundani,rank$\}$@bv.tum.de}}
 
  \maketitle

\begin{abstract}
The finite cell method (FCM) belongs to the class of immersed boundary methods, and combines the fictitious domain approach with high-order approximation, adaptive integration and weak imposition of unfitted Dirichlet boundary conditions. For the analysis of complex geometries, it circumvents expensive and potentially error-prone meshing procedures, while maintaining high rates of convergence. The present contribution provides an overview of recent accomplishments in the FCM with applications in structural mechanics. First, we review the basic components of the technology using the $p$- and B-spline versions of the FCM. Second, we illustrate the typical solution behavior for linear elasticity in 1D. Third, we show that it is straightforward to extend the FCM to nonlinear elasticity. We also outline that the FCM can be extended to applications beyond structural mechanics, such as transport processes in porous media. Finally, we demonstrate the benefits of the FCM with two application examples, i.e.~the vibration analysis of a ship propeller described by T-spline CAD surfaces and the nonlinear compression test of a CT-based metal foam.
\end{abstract}

\pagestyle{empty}
\thispagestyle{fancy}
\lhead{}
\chead{}
\rhead{}
\lfoot{\scriptsize This is a pre-print of an article published in Bader~M., Bungartz~HJ., Weinzierl~T.\ (eds) Advanced Computing. Lecture Notes in Computational Science and Engineering, vol 93, 2013. The final authenticated version is available online at: https://doi.org/10.1007/978-3-642-38762-3\_1}
\cfoot{}
\rfoot{}

\flushbottom

\section{Introduction}

Structural analysis with standard finite elements requires the discretization of the domain of interest into a finite element mesh, whose boundaries conform to the physical boundaries of the structure \cite{Bathe:96.1,Hughes:00.1}. While this constraint can be easily achieved for many applications in structural mechanics, it constitutes a severe bottleneck when highly complex geometries are concerned. An alternative pathway that avoids time-consuming mesh generation is provided by embedded domain methods, also known as immersed boundary methods \cite{Neittaanmaeki:95.1,Peskin:02.1,Mittal:05.1,Loehner:08.1}. Their main idea consists of the extension of the physical domain of interest $\Omega_{\textit{phys}}$ beyond its potentially complex boundaries into a larger embedding domain of simple geometry $\Omega$, which can be meshed easily by a structured grid (see Fig.~1). The finite cell method (FCM) \cite{Parvizian:07.1,Duester:08.1,Schillinger:12.4} is an embedded domain method, which combines the fictitious domain approach \cite{Bishop:03.1,Glowinski:07.1,Ramiere:07.1} with higher-order basis functions \cite{Szabo:04.1,Schillinger:12.2}, adaptive integration and weak imposition of unfitted Dirichlet boundary conditions \cite{Zhu:98.1,Fernandez:04.1,Ruess:13.1}. To preserve consistency with the original problem, the influence of the fictitious domain extension $\Omega_{\textit{fict}}$ is extinguished by penalizing its material parameters. For smooth problems of linear elasticity, the FCM has been shown to maintain exponential rates of convergence in the energy norm and thus allows for accurate structural analysis irrespective of the geometric complexity involved \cite{Rank:09.1}. Moreover, it can be well combined with image-based geometric models typical for applications from biomechanics and material science \cite{Duester:08.1,Ruess:12.1,Schillinger:12.1}. Within the framework of the FCM for structural analysis, the following aspects have been examined so far: Topology optimization \cite{Parvizian:12.1}, thin-walled structures \cite{Rank:12.1}, local refinement strategies \cite{Schillinger:12.1,Schillinger:10.1,Schillinger:11.1,Schillinger:12.3}, weak boundary conditions \cite{Schillinger:12.2,Ruess:12.1}, homogenization of porous and cellular materials \cite{Duester:12.1}, geometrically nonlinear problems \cite{Schillinger:11.2,Schillinger:12.2}, and computational steering \cite{Knezevic:11.1,Knezevic:11.2,Yang:12.1,Yang:12.2}.

The present contribution provides an overview of recent accomplishments in the finite cell method for structural mechanics. It is organized as follows: \textbf{Section~2} provides a short introduction to the basic components of the finite cell method. \textbf{Section~3} outlines the typical solution behavior for linear elastic problems and highlights important numerical properties. \textbf{Section~4} shows the extension of the finite cell method to nonlinear elasticity. \textbf{Section~5} outlines that FCM can be extended to problems beyond structural mechanics by the example of transport processes in porous media. \textbf{Section~6} presents two application oriented numerical examples in three-dimensions, based on CAD and image-based geometric models. In \textbf{Section~7}, we conclude our presentation by a short summary and an outlook to future research.

\begin{figure*}[t!]
\centering
\psfrag{O1}{\footnotesize{$\Omega_{\textit{phys}}$}}
\psfrag{t}{\footnotesize{$\boldsymbol{t}$}}
\psfrag{dO}{\footnotesize{$\partial\Omega$}}
\psfrag{GN}{\footnotesize{$\Gamma_N$}}
\psfrag{GD}{\footnotesize{$\Gamma_D$}}
\psfrag{O2}{\hspace{-0.008\textwidth}\footnotesize{$\Omega_{\textit{fict}}$}}
\psfrag{O1+O2}{\footnotesize{$\Omega$}=\footnotesize{$\Omega_{\textit{phys}}$}+\footnotesize{$\Omega_{\textit{fict}}$}}
\psfrag{a1}{\footnotesize{$\alpha=1.0$}}
\psfrag{a2}{\footnotesize{$\alpha \ll 1.0$}}
\includegraphics[width=0.98\textwidth]{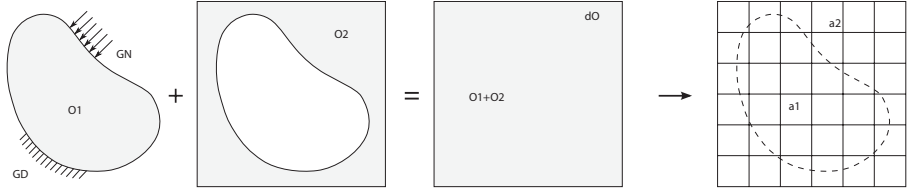}
\caption{The fictitious domain concept: The physical domain $\Omega_{\textit{phys}}$ is extended by the fictitious domain $\Omega_{\textit{fict}}$ into an embedding domain $\Omega$ to allow easy meshing of complex geometries. The influence of $\Omega_{\textit{fict}}$ is penalized by $\alpha$.} 
\end{figure*}

\section{A Brief Review of the Finite Cell Method}

The following review provides a brief introduction to the main components, i.e.~the fictitious domain concept, a higher-order approximation basis, adaptive integration and unfitted Dirichlet boundary conditions. We follow the presentation given in \cite{Schillinger:12.2}, focusing on the $p$- and the B-spline versions of the FCM.

\subsection{The Fictitious Domain Concept}

In the finite cell method, the domain to be analyzed is called the embedding domain $\Omega$, which consists of the physical domain of interest $\Omega_{\textit{phys}}$ and the fictitious domain extension $\Omega_{\textit{fict}}$ as shown in Fig.~1. Analogous to standard finite element methods (FEM), the finite cell method for linear elastic problems is derived from the principle of virtual work
\begin{align}
\delta W\left(\boldsymbol{u},\delta \boldsymbol{u}\right) = \int_{\Omega} \boldsymbol{\sigma} \; : \; (\nabla_{\hspace{-0.069cm}sym} \, \delta \boldsymbol{u}) \;\; dV - \int_{\Omega_{\textit{phys}}} \delta \boldsymbol{u} \cdot \boldsymbol{b}  \; dV - \int_{\Gamma_N} \delta  \boldsymbol{u} \cdot \boldsymbol{t} \; dA \;\; = \;\;  0
\end{align}

\noindent where $\boldsymbol{\sigma}$, $\boldsymbol{b}$, $\boldsymbol{u}$, $\delta \boldsymbol{u}$ and $\nabla_{\hspace{-0.069cm}sym}$ denote the Cauchy stress tensor, body forces, displacement vector, test function and the symmetric part of the gradient, respectively \cite{Bathe:96.1,Hughes:00.1}. Neumann boundary conditions are specified over the boundary of the embedding domain $\partial\Omega$, where tractions are zero by definition, and over $\Gamma_N$ of the physical domain by traction vector $\boldsymbol{t}$ (see Fig.~1). The elasticity tensor $\boldsymbol{C}$ \cite{Bathe:96.1,Hughes:00.1} relating stresses and strains
\begin{equation}
\boldsymbol{\sigma} = \alpha \boldsymbol{C} \, : \, \boldsymbol{\varepsilon}
\end{equation}

\noindent is complemented by a scalar factor $\alpha$, which reads
\begin{equation}
\alpha\left(\boldsymbol{x}\right) = \begin{cases}
  1.0    & \forall \boldsymbol{x} \in \Omega_{\textit{phys}} \\
  10^{-q}    & \forall \boldsymbol{x} \in \Omega_{\textit{fict}}
\end{cases}
\end{equation}

\noindent penalizing the contribution of the fictitious domain. In $\Omega_{\textit{fict}}$, $\alpha$ must be chosen as small as possible, but large enough to prevent extreme ill-conditioning of the stiffness matrix \cite{Parvizian:07.1,Duester:08.1}. Typical values of $\alpha$ range between $10^{-4}$ and $10^{-15}$. 

Using a structured grid of high-order elements (see Fig.~1), which will be called finite cells in the following, kinematic quantities are discretized as
\begin{equation}
\boldsymbol{u}= \sum_{a=1}^n N_a \boldsymbol{u}_a
\end{equation}
\begin{equation}
\delta\boldsymbol{u}= \sum_{a=1}^n N_a \delta \boldsymbol{u}_a
\end{equation}

\noindent The sum of $N_a$ denotes a finite set of $n$ higher-order shape functions, and $\boldsymbol{u}_a$ and $\delta\boldsymbol{u}_a$ the corresponding vectors of unknown coefficients. Following the standard Galerkin approach \cite{Bathe:96.1,Hughes:00.1}, inserting (4) and (5) into the weak form (1) produces a discrete finite cell representation
\begin{equation}
\boldsymbol{K} \boldsymbol{u} = \boldsymbol{f}
\end{equation}

\noindent with stiffness matrix $\boldsymbol{K}$ and load vector $\boldsymbol{f}$. Due to the similarity to standard FEM, the implementation of FCM can make use of existing techniques.

\begin{figure*}[t!]
\centering
\psfrag{xi}{\footnotesize{$\xi$}}
\psfrag{N1}{\footnotesize{$N_1$}}
\psfrag{N2}{\footnotesize{$N_2$}}
\psfrag{p2}{\footnotesize{$\phi_2$}}
\psfrag{p3}{\footnotesize{$\phi_3$}}
\psfrag{p4}{\footnotesize{$\phi_4$}}
\psfrag{p5}{\footnotesize{$\phi_5$}}
\includegraphics[width=0.32\textwidth]{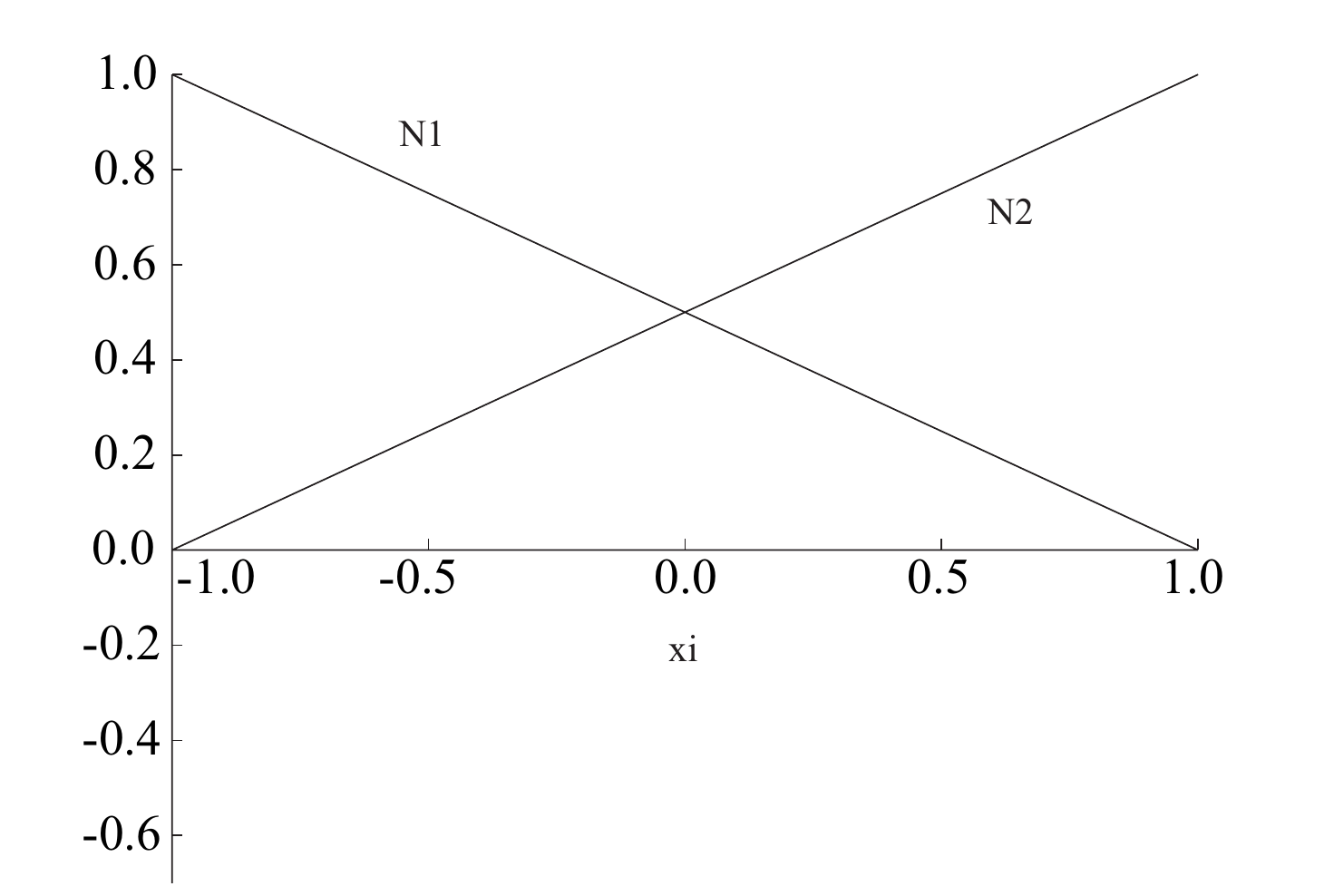}
\includegraphics[width=0.32\textwidth]{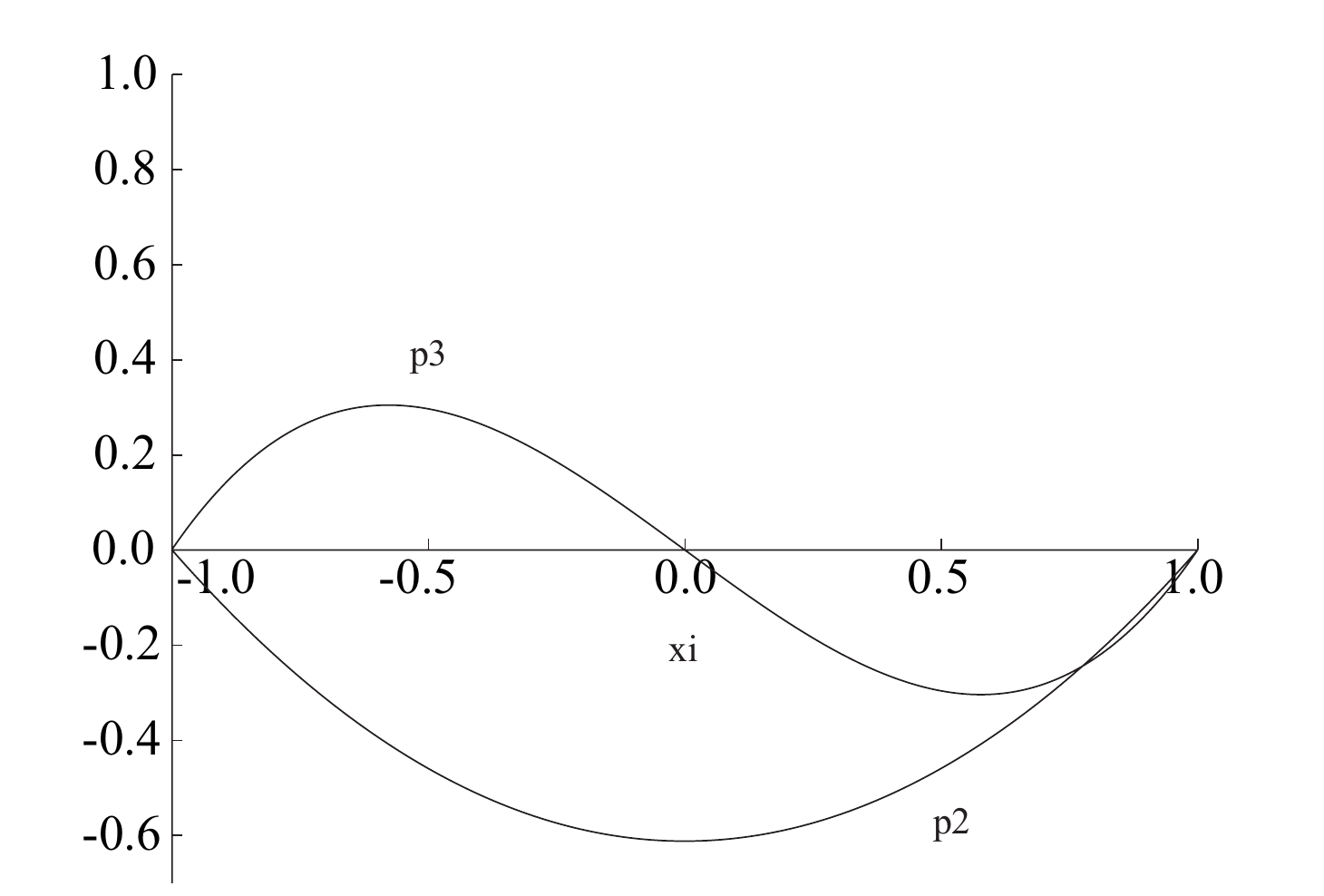}
\includegraphics[width=0.32\textwidth]{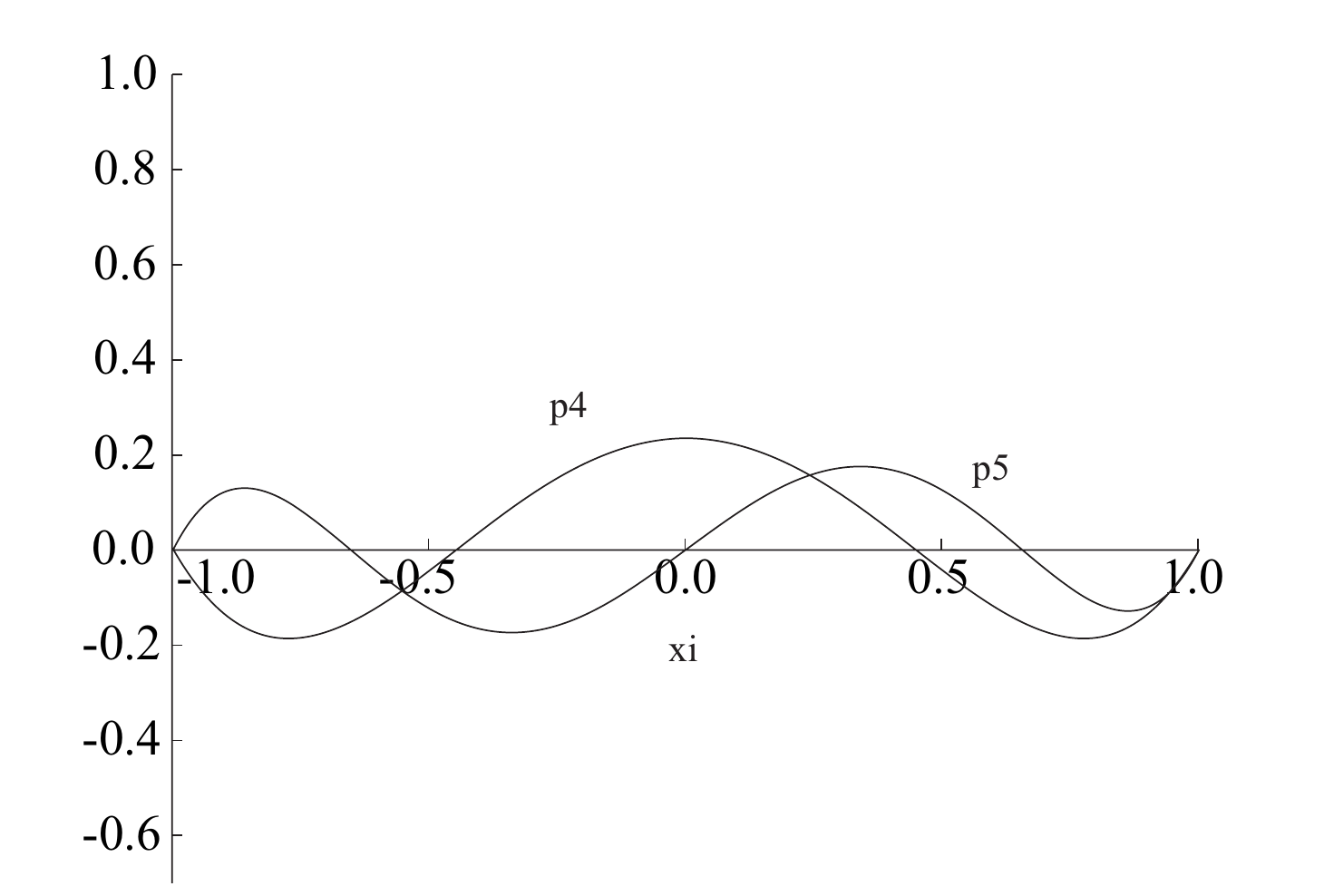}
\caption{Linear nodal modes $N_j$, $j$=$1,2$ and the first 4 integrated Legendre basis functions $\phi_j$, $j$=$2,...,5$ of the 1D $p$-version basis in the parameter space $\xi$.}
\end{figure*}
\begin{figure*}[t!]
\psfrag{0}{\scriptsize{0}}
\psfrag{1}{\scriptsize{1}}
\psfrag{2}{\scriptsize{2}}
\psfrag{3}{\scriptsize{3}}
\psfrag{4}{\scriptsize{4}}
\psfrag{a1}{\scriptsize{$\xi$}}
\psfrag{a2}{\scriptsize{$\eta$}}
\parbox{1.0\textwidth}{\centering\includegraphics[width=0.36\textwidth]{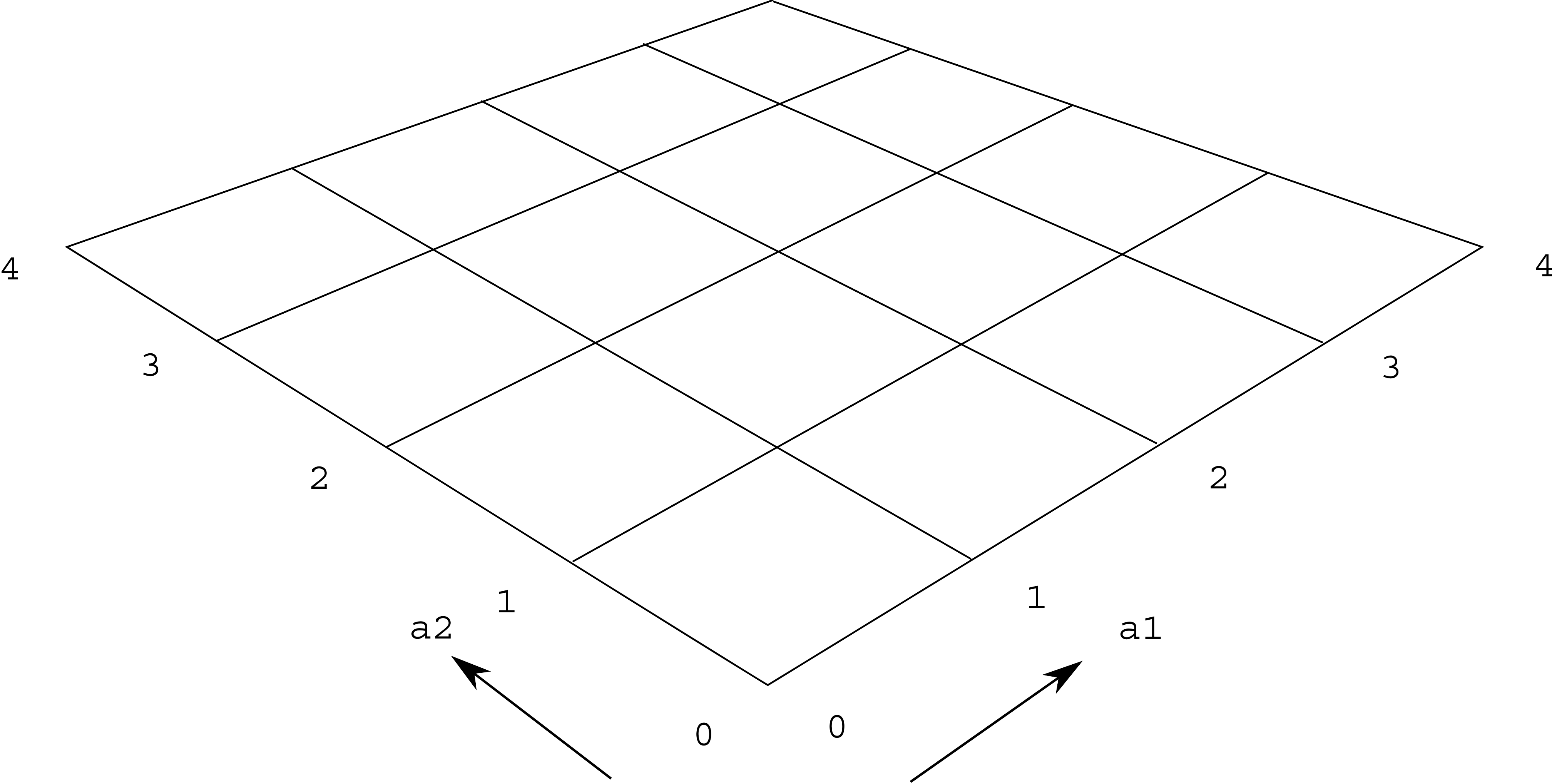}
\hspace{1.cm} \includegraphics[width=0.38\textwidth]{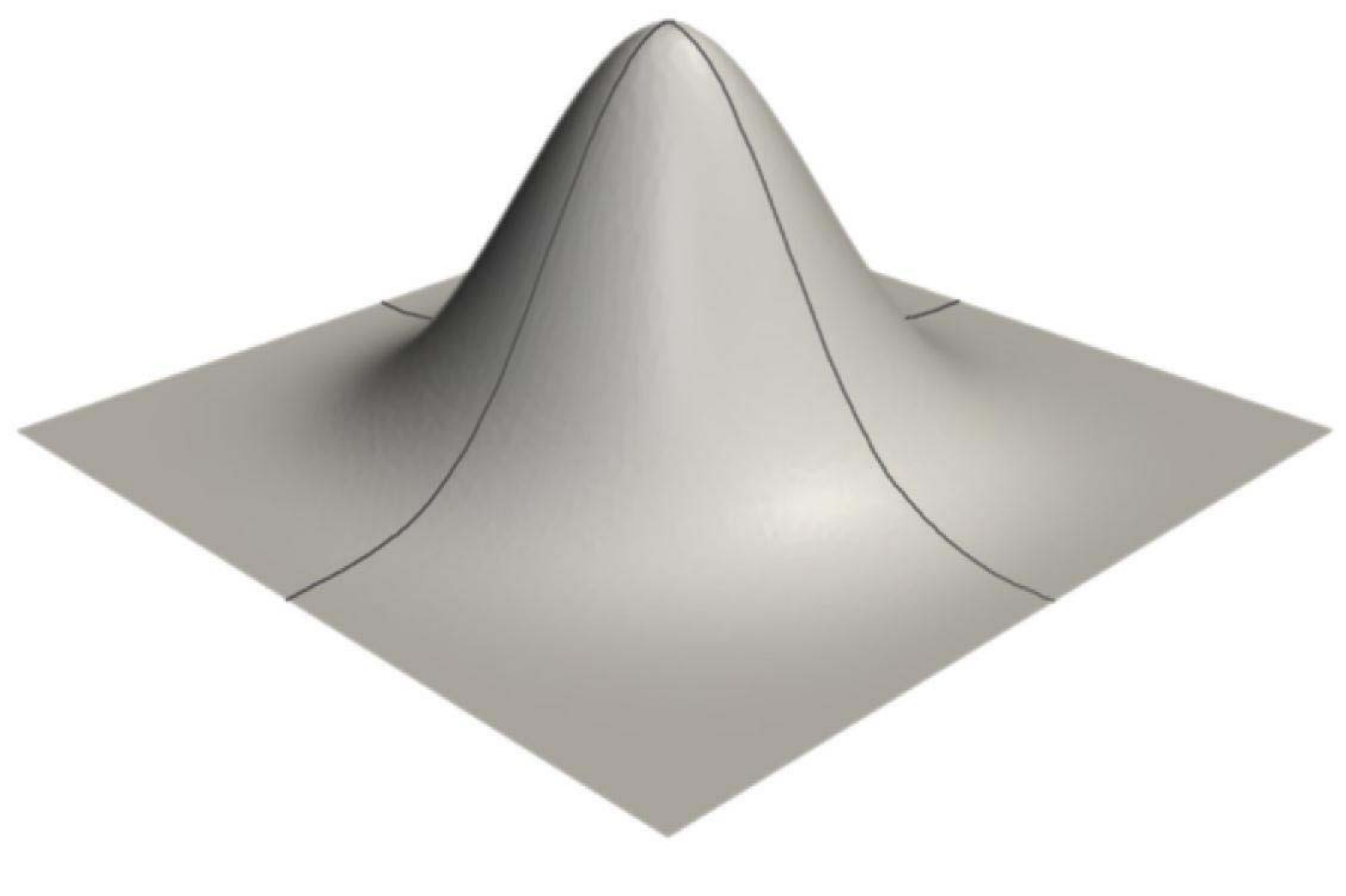} }
\caption{Knot span cells in the parameter space $\{\xi,\eta\}$ (left) and corresponding bi-variate cubic B-spline (right).}
\end{figure*}

\subsection{Higher-order Approximation of Solution Fields}

The high-order basis originally applied in the FCM \cite{Parvizian:07.1,Duester:08.1} uses a regular mesh of elements of the \textit{p-version of the FEM}, whose formulation is based on $C^0$ integrated Legendre polynomials \cite{Szabo:04.1,Szabo:91.1}. Corresponding basis functions in 1D are plotted in Fig.~2. The basis is hierarchical, so that an increase of the polynomial degree $p$ by 1 is achieved by the addition of another function $\phi_j$. Corresponding higher-dimensional bases can be constructed by tensor products of the 1D case. To limit the number of additional unknowns in 2D and 3D, the so-called \textit{trunk space} is used instead of the full tensor product basis \cite{Szabo:04.1,Szabo:91.1}.

The \textit{B-spline version of the FCM} has been recently established as a suitable alternative \cite{Schillinger:12.2,Schillinger:10.1,Schillinger:11.1,Schillinger:12.3}. Its formulation is based on higher-order and smooth B-spline basis functions \cite{Piegl:97.1,Rogers:01.1}, whose numerical advantages have been recently demonstrated in the context of isogeometric analysis \cite{Hughes:05.1,Cottrell:09.1}. We use a single uniform B-spline patch, whose basis functions consist of uniform B-splines constructed from equidistant knots \cite{Piegl:97.1,Rogers:01.1} and can be interpreted as translated copies of each other \cite{Hoellig:03.1}. Corresponding multivariate B-spline basis functions are obtained by taking the tensor product of the univariate components in each parametric direction. An example of a two-dimensional knot span structure and a corresponding bi-cubic uniform B-spline are shown in \mbox{Fig.~3}. Each knot span can be identified as a quadrilateral or hexahedral finite cell, respectively, with full Gaussian integration \cite{Schillinger:10.1,Schillinger:11.1}. The physical coordinates of the FCM grid can be generated from a simple linear transformation of the parametric space \cite{Schillinger:12.2}. 

\begin{figure*}[t!]
\centering
\psfrag{Initial}{\footnotesize{Finite cell mesh}}
\psfrag{k0}{\footnotesize{\textit{k=0}}}
\psfrag{k1}{\footnotesize{\textit{k=1}}}
\psfrag{k2}{\footnotesize{\textit{k=2}}}
\psfrag{k3}{\footnotesize{\textit{k=3}}}
\psfrag{k4}{\footnotesize{\textit{k=4}}}
\psfrag{k5}{\footnotesize{\textit{k=5}}}
\psfrag{G1}{\footnotesize{with geometric}}
\psfrag{G2}{\footnotesize{boundary}}
\includegraphics[width=0.92\textwidth]{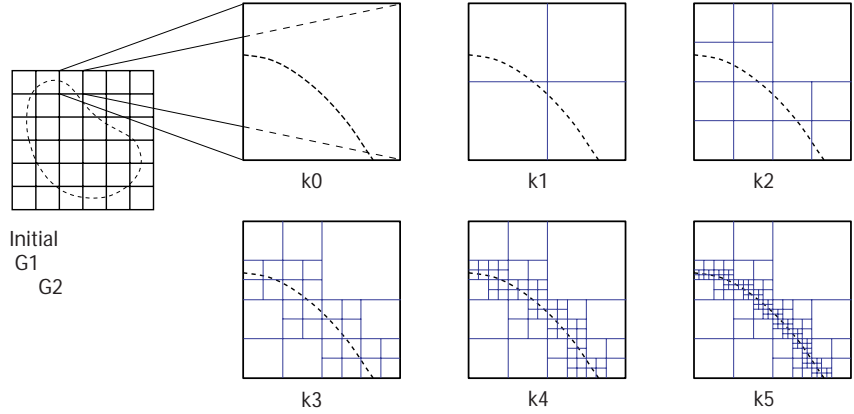}
\caption{2D sub-cell structure (thin blue lines) for adaptive integration of finite cells (bold black lines) that are cut by the geometric boundary (dashed line).} 
\end{figure*}

\subsection{Adaptive Integration}

The accuracy of numerical integration by Gauss quadrature \cite{Bathe:96.1,Hughes:00.1} is considerably influenced by discontinuities within cells introduced by the penalization parameter $\alpha$ of (3) \cite{Parvizian:07.1,Duester:08.1}. Therefore, the FCM uses composed Gauss quadrature in cells cut by geometric boundaries, based on a hierarchical decomposition of the original cells \cite{Duester:08.1,Schillinger:12.2,Schillinger:12.1}. In two dimensions, the sub-cell structure can be built up in the sense of a quadtree (see \mbox{Fig.~4}) \cite{Samet:06.1}. Starting from the original finite cell of level \textit{k=0}, each sub-cell of level \textit{k=i} is first checked whether it is cut by a geometric boundary. If true, it is replaced by 4 equally spaced cells of level \textit{k=i+1}, each of which is equipped with ($p$+1)$\times$($p$+1) Gauss points. Partitioning is repeated for all cells of current level $k$, until a predefined maximum depth \textit{k=m} is reached. The quadtree approach can be easily adjusted to 1D or 3D by binary trees or octrees, respectively. It is easy to implement and keeps the regular grid structure of the FCM. To clearly distinguish between finite cell and sub-cell meshes, finite cells are plotted in black and integration sub-cells are plotted in blue lines throughout this paper (see Fig.~4).

\subsection{Imposition of Unfitted Boundary Conditions}

For complex domains, boundary conditions are defined along geometric boundaries cutting arbitrarily through finite cells. Neumann boundary conditions can be incorporated by simple integration over the Neumann boundary $\Gamma_N$ (see (1)). Dirichlet boundary conditions require an imposition in a weak sense by variational techniques such as the penalty method \cite{Zhu:98.1,Babuska:72.1}, the Lagrange multiplier method \cite{Glowinski:07.1,Fernandez:04.1,Zienkiewicz:00.1} or Nitsche's method \cite{Hansbo:02.1,Bazilevs:07.2,Embar:10.1}. In the FCM, Nitsche's method is usually preferred \cite{Schillinger:12.2,Ruess:12.1,Ruess:13.1}, since it does not introduce additional unknowns, leads to a symmetric, positive definite stiffness matrix and satisfies variational consistency in the sense that solutions of the weak form can be shown to be solutions of the original boundary value problem.

From a practical point of view, the integration over unfitted boundaries is accomplished by introducing a triangular mesh of the boundary surfaces. Generating a triangulation of a 3D surface is a standard task, for which a variety of efficient algorithms and tools are available. In particular, it is orders of magnitude less complex and less expensive than the generation of a full volumetric discretization of a complex 3D object. Corresponding mesh generation in the framework of the FCM for a boundary representation of solids and for voxel-based data obtained from CT scans is addressed in detail in \cite{Duester:08.1}.

\begin{figure*}[t!]
\centering
\psfrag{O1}{\footnotesize{$\Omega_{\textit{phys}}$}}
\psfrag{O2}{\footnotesize{$\Omega_{\textit{fict}}$}}
\psfrag{L}{\footnotesize{$L$}}
\psfrag{L2}{\footnotesize{$4/3\,L$}}
\psfrag{L3}{\footnotesize{$2/3\,L$}}
\psfrag{x}{\footnotesize{$X$}}
\psfrag{du}{\footnotesize{$\Delta u$}}
\psfrag{f}{\footnotesize{$f_{sin}$}}
\psfrag{a1}{\footnotesize{Discretization with 2 $p$-version cells:}}
\psfrag{b1}{\footnotesize{Discretization with 11 B-spline knot span cells:}}
\psfrag{a2}{\parbox{5cm}{\footnotesize{Parameters:}\\[0.26cm]
												 \footnotesize{Young's modulus \textit{E=1.0}}\\[0.046cm]
                         \footnotesize{Poisson's ratio $\nu=0.0$}\\[0.046cm]
                         \footnotesize{Penalization parameter $\alpha=10^{-8}$}\\[0.046cm]
                         \footnotesize{Area \textit{A=1.0}}\\[0.046cm]
                         \footnotesize{Length \textit{L=1.0}}\\[0.046cm]
                         \footnotesize{Displacement load $\Delta u$=$0.02$}\\[0.046cm]
                         \footnotesize{Sine load $f_{\text{sin}}=1/20 \sin{(4\pi X)}$}}}
\includegraphics[width=0.9\textwidth]{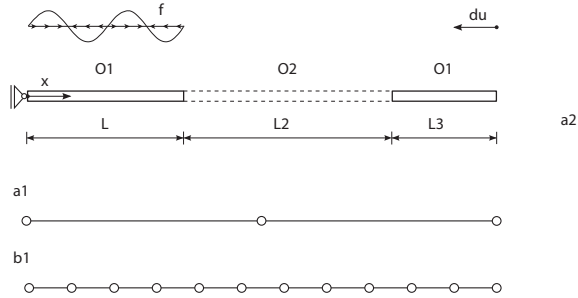}
\caption{Uni-axial rod with geometric boundaries at $X=1.0$ and $X=2 \frac{1}{3}$.}
\end{figure*}

\section{Basic Numerical Properties of the FCM}

For the illustration of the typical solution behavior, a linear elastic uni-axial rod is examined, for which geometry, material and boundary conditions are specified in Fig.~5. Its middle part represents the fictitious domain $\Omega_{\textit{fict}}$, whose Young's modulus $E$ is penalized with parameter $\alpha$\textit{=10$^{-8}$}. The example approximates the situation of two separate rods. The right one undergoes a rigid body movement $\Delta u$ and the left one is subjected to a sine load $f_{\text{sin}}$. The FCM discretizations considered consist of 2 $p$-version finite cells and 11 knot span cells as shown in Fig.~5. Due to the different construction of the bases, the B-spline version requires a denser knot span grid than the $p$-version in order to achieve a comparable amount of degrees of freedom (dofs). For all computations of this section, adaptive sub-cells of depth $k$=20 are used to minimize the integration error in cells cut by geometric boundaries. 

\subsection{Smooth Extension of Solution Fields}

The $p$- and B-spline versions of the FCM produce solution fields, which extend smoothly into the fictitious domain despite the discontinuities of the analytical solution. This is illustrated in Fig.~6, which compares the analytical strains to the numerical strains of the $p$- and B-spline versions. The importance of the smooth extension of the FCM solution for the overall convergence behavior of the finite cell method can be explained with the help of the penalty parameter $\alpha$ in conjunction with the best approximation property to the total strain energy \cite{Strang:73.1} (see for example \cite{Schillinger:12.2,Schillinger:12.3} for details). It should be noted that the large difference between analytical and FCM solution fields in the fictitious domain (see Fig.~6) is completely irrelevant, since we are only interested in the physical solution.

\begin{figure*}[t!]
\centering
\includegraphics[width=0.62\textwidth]{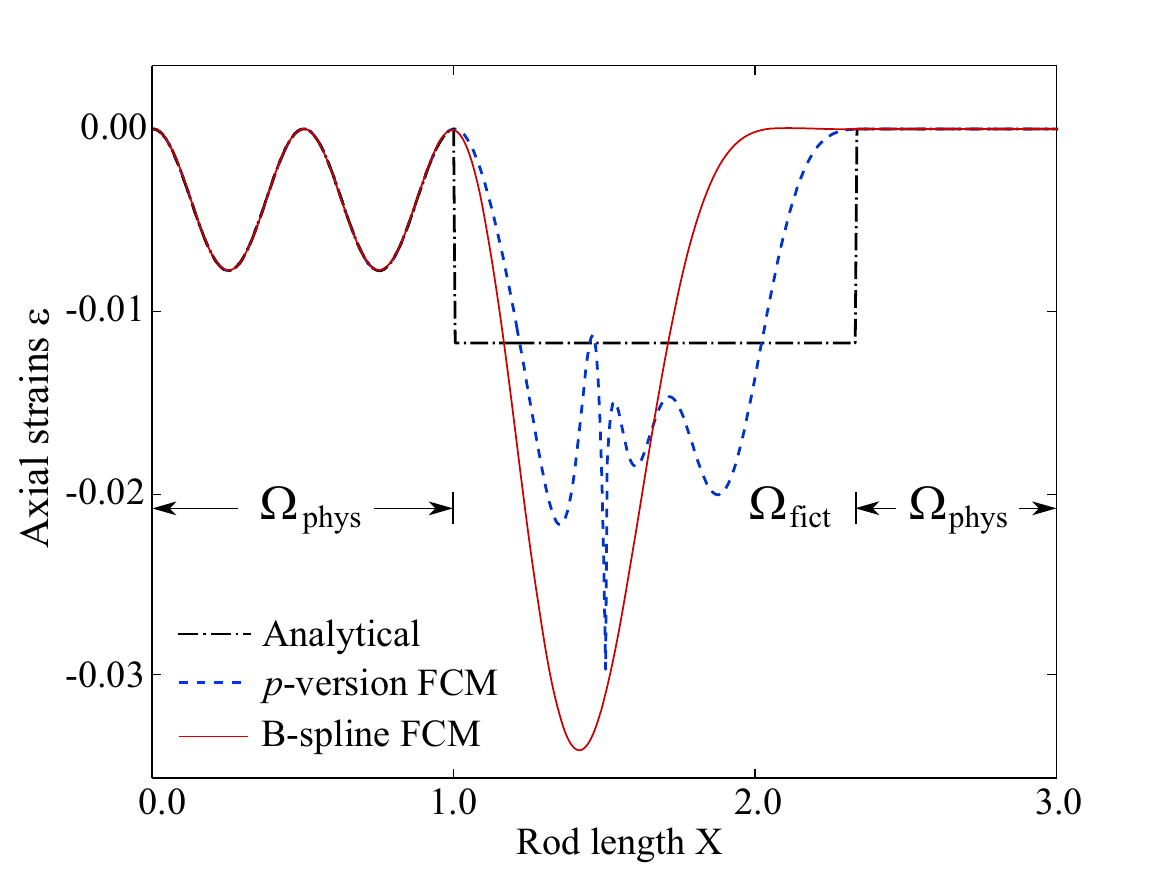}
\caption{Smooth extension of the FCM solutions ($p$=15) \mbox{vs.} discontinuous analytical solution ($\alpha$=$10^{-8}$) for the linear elastic strains of the 1D example.}
\end{figure*}

\subsection{Exponential Rates of Convergence}

Figure~7 shows the convergence behavior of the presented FCM schemes, if the polynomial degree of the discretizations given in Fig.~5 is increased from $p$=1 to 15. Both the $p$-version and the B-spline version of the FCM converge exponentially. The penalization parameter \textit{$\alpha$=$10^{-8}$} as well as the integration error in cut cells lead to a flattening of the convergence curves. The present example shows that the $p$- and B-spline bases exhibit an equivalent solution behavior within the FCM and achieve a comparable performance in terms of error level, rates of convergence and flattening of the convergence curve, although their high-order approximation bases are very different. Further numerical benchmarks in higher dimensions can be found in \cite{Schillinger:12.2,Schillinger:12.3} that show optimal rates of convergence under $h$-refinement, the stability and accuracy of weak boundary conditions and the competitive quality of the solution and its derivatives along the geometric boundaries in cut cells with respect to standard body-fitted finite element methods.

\begin{figure*}[t!]
\centering
\includegraphics[width=0.62\textwidth]{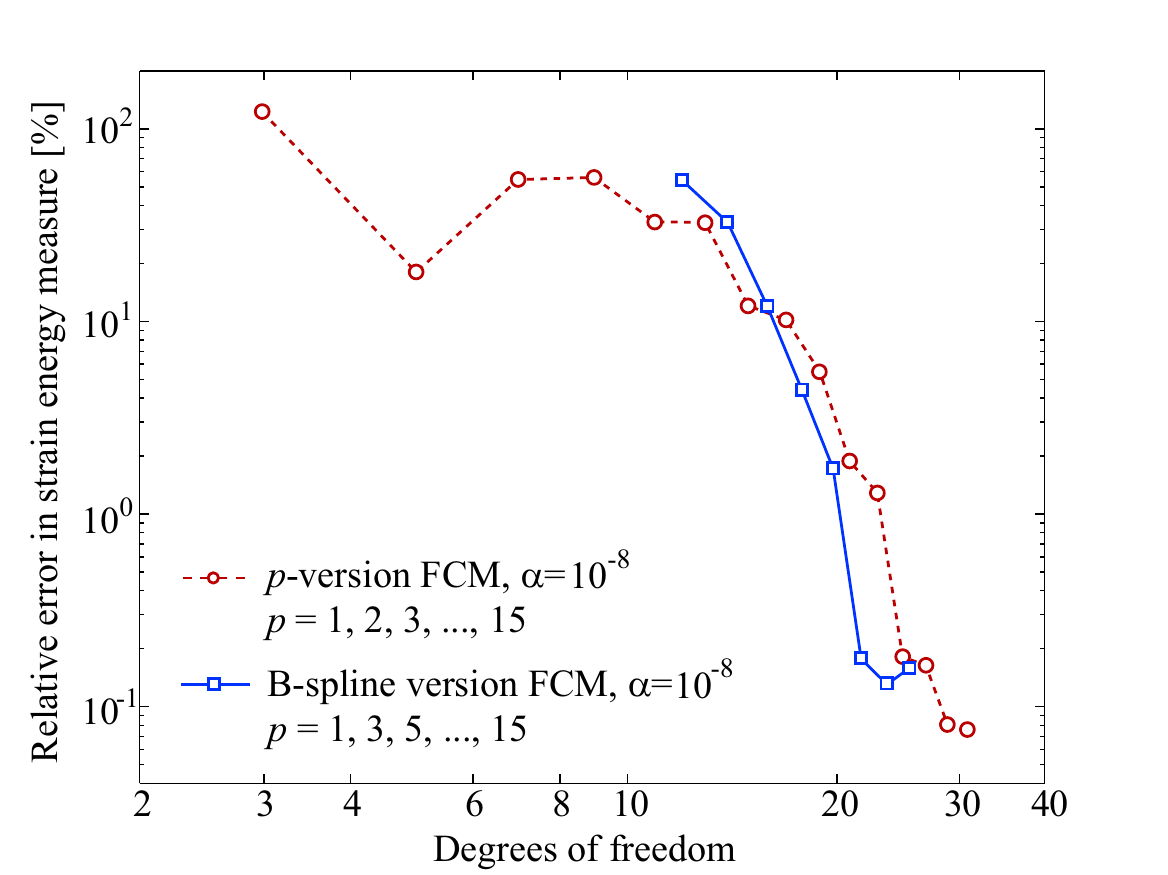}
\caption{Convergence in strain energy obtained by $p$-refinement in the given $p$-version and B-spline discretizations of the linear elastic 1D example.}
\end{figure*}

\section{Extension to Nonlinear Elasticity}

The finite cell method can be extended to geometrically nonlinear elasticity on the basis of the logarithmic strain measure \cite{Bonet:08.1} and the Hencky hyperelastic material model \cite{deSouzaNeto:08.1}. An extensive review of the mathematical model and the pertinent continuum mechanics in the framework of the FCM can be found in \cite{Schillinger:12.2}. In the present scope, we focus on our 1D rod example of Fig.~5, for which the corresponding geometrically nonlinear formulation simplifies to 
\begin{equation}
\Psi = \alpha \, \frac{E}{2} \left(\ln \lambda \right)^2
\end{equation}
\begin{equation}
\sigma = \alpha \, \frac{E}{J} \ln \lambda
\end{equation}
\begin{equation}
c = \alpha \, \frac{E}{J} - 2 \sigma
\end{equation}

\noindent with axial stretch $\lambda$, the determinant of the deformation gradient $J=\lambda^{1-2\nu}$ and the strain energy function $\Psi$ \cite{Bonet:08.1}. To illustrate the influence of large deformations within the fictitious domain, the sine load $f_{\text{sin}}$ of Fig.~5 is neglected for now and the prescribed displacement is set to a large value of $\Delta u$=1.0. The physical stresses should be zero, since a rigid body movement of the right part of the rod is approximated. The exact stress solution, which can be derived analytically according to \cite{Schillinger:10.1}, is plotted in Fig.~8a for 10 displacement load increments between 0 and $\Delta u$ and $\alpha=10^{-5}$. In the following, all numerical examples are computed with the B-spline version of the FCM, but equivalent results can be derived for the $p$-version (see \cite{Schillinger:12.2}).

\begin{figure*}[h!]
\centering
\subfloat[Analytical stress solution.] {\hspace{2cm}
\includegraphics[width=0.52\textwidth]{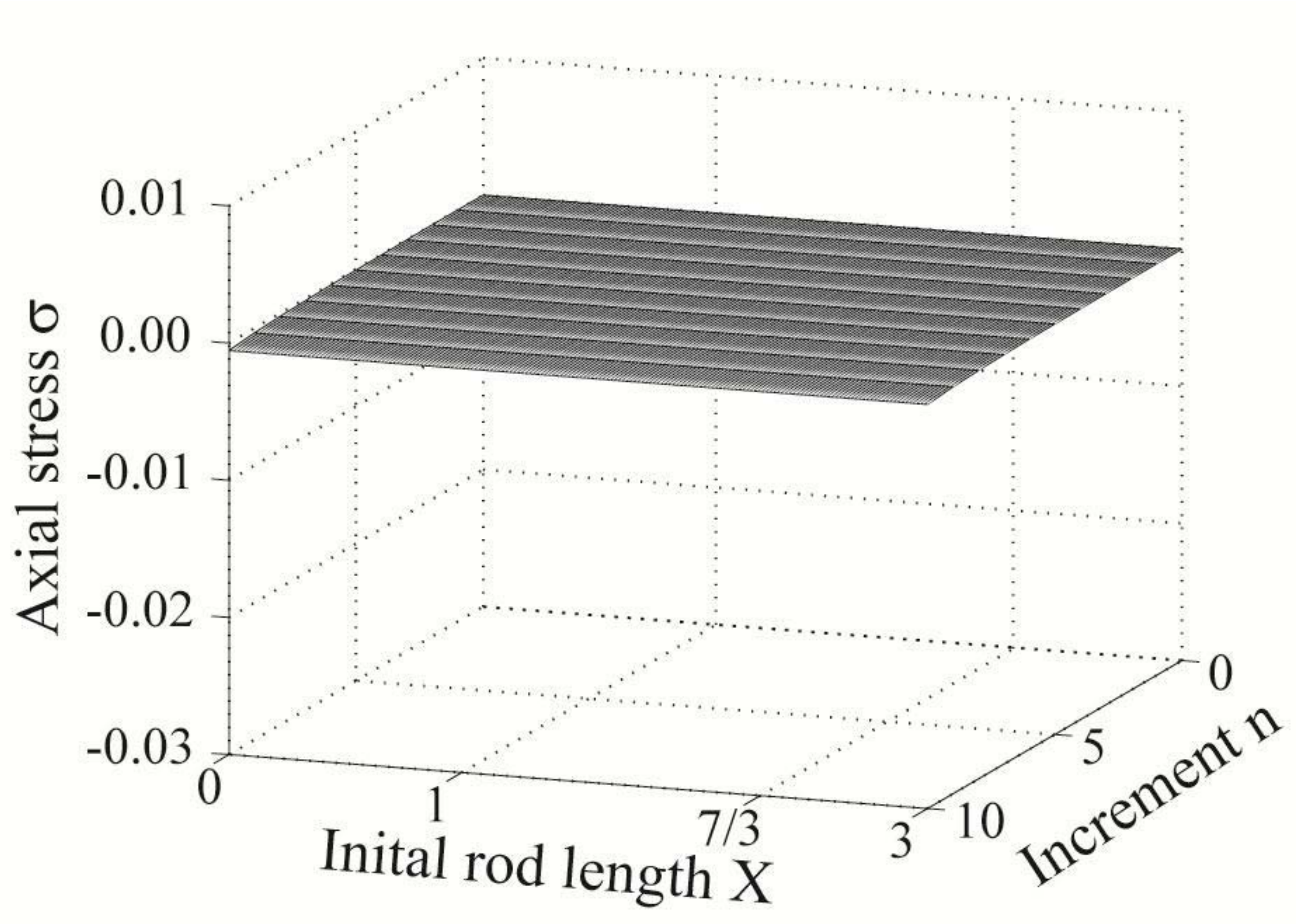} \hspace{2cm}}\\[0.42cm]

\subfloat[Computed with the standard B-spline version of the FCM.] {\hspace{2cm}
\includegraphics[width=0.52\textwidth]{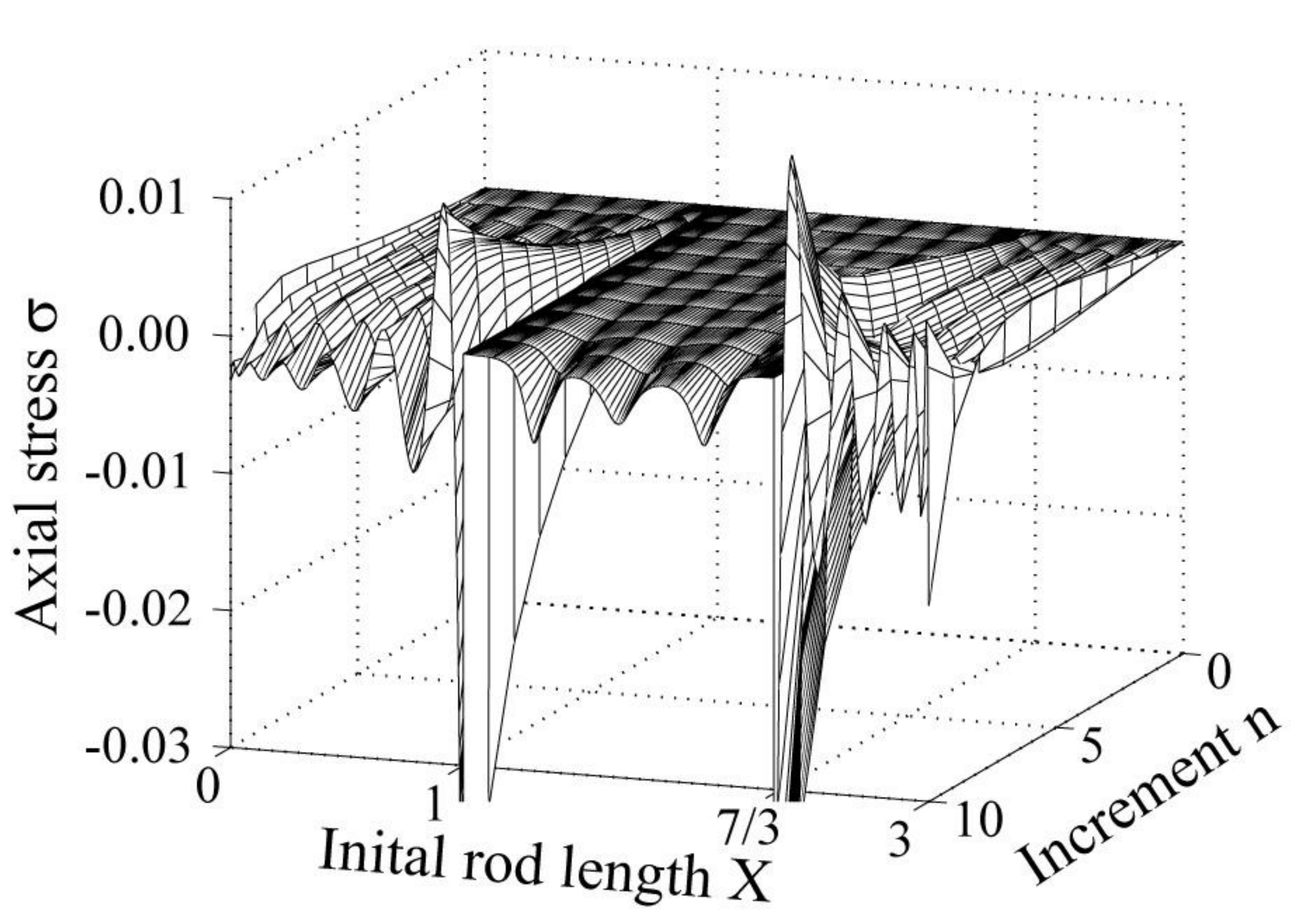}\hspace{2cm}}\\[0.42cm]

\subfloat[Computed with the B-spline version of the FCM and deformation resetting.] {\hspace{2.8cm}
\includegraphics[width=0.52\textwidth]{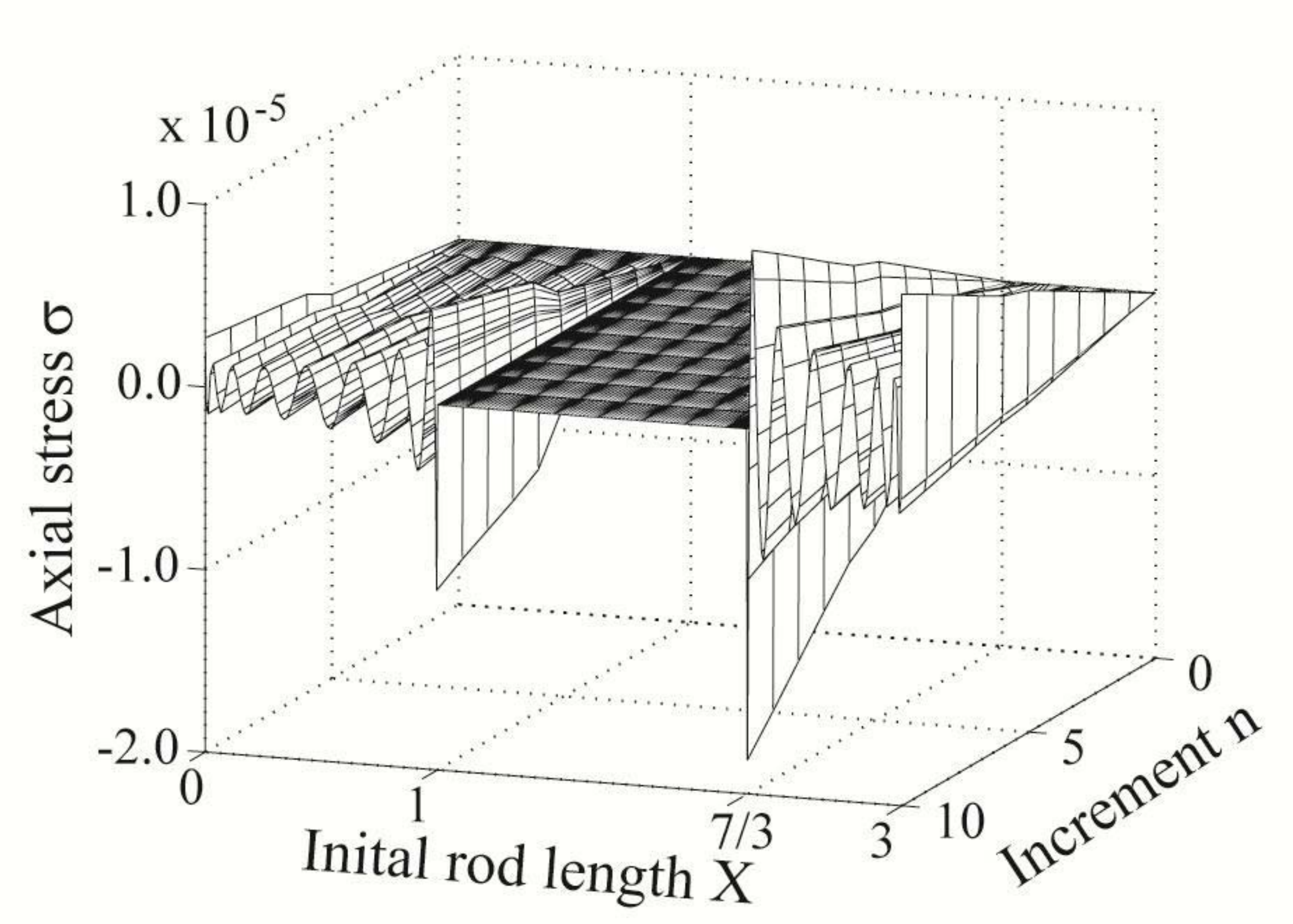}\hspace{2.8cm} }
\caption{The stress solution of the geometrically nonlinear rod is obtained without and with deformation resetting (16 knot span elements of $p$=15). Note that deformation resetting reduces the stress oscillations by three orders of magnitude.}
\end{figure*}

\subsection{The Standard FCM Formulation}

The standard FCM is based on the application of the same geometrically nonlinear formulation over the complete embedding domain $\Omega$. However, numerical experiments reveal that for $\alpha$ smaller than $10^{-5}$, the determinant of the deformation gradient falls below zero at some integration point within $\Omega_{\textit{fict}}$, which inevitably terminates the computation. With $\alpha$ as large as $10^{-5}$, the penalization of (3) is unable to sufficiently eliminate the influence of $\Omega_{\textit{fict}}$, so that a considerable modeling error is introduced. In addition, nonlinear strains increasingly outweigh the penalization by $\alpha$, since they are able to grow without bounds \cite{Schillinger:12.2}. The corresponding stress solution obtained with 16 knot span cells in the sense of Fig.~5 is plotted in Fig.~8b. It exhibits large oscillations throughout the discontinuous cells and the corresponding convergence deteriorates to a low algebraic rate (see Fig.~9). The standard FCM formulation thus suffers from a conflict of interest between stable analysis (increase of $\alpha$) on the one hand and a reduction of the contribution of $\Omega_{\textit{fict}}$ (decrease of $\alpha$) on the other.

\subsection{A Modified Formulation based on Deformation Resetting}

Numerical experiments reveal that problems with the uniqueness of the deformation map occur at the location of maximum deformation within the fictitious domain $\Omega_{\textit{fict}}$. This motivates the following simple manipulation after each Newton iteration $i$
\begin{equation}
\boldsymbol{\varphi}^i({\bold X}) = \begin{cases}
\;\boldsymbol{x}^i   & \text{deformed configuration} \;\; \forall \boldsymbol{X} \in \Omega_{\textit{phys}} \\[0.1cm]
\;\boldsymbol{X} & 
\begin{array}{l}
\text{reset to reference} \\[-0.1cm]
\text{configuration} \;\;  \forall \boldsymbol{X} \in \Omega_{\textit{fict}} 
\end{array}
\end{cases}
\end{equation}

\begin{figure}[t!]
\centering
\includegraphics[width=0.62\textwidth]{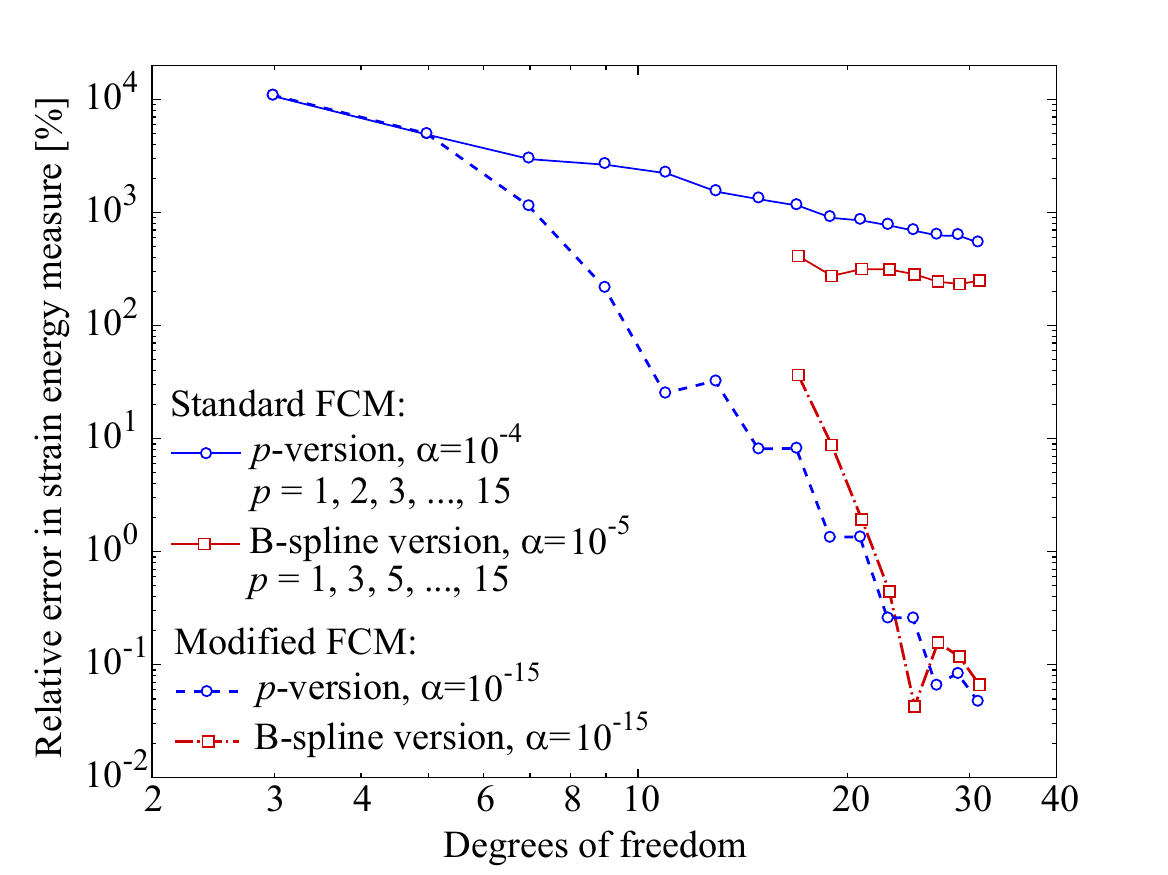}
\caption{Convergence without and with deformation resetting for the nonlinear rod.}
\end{figure}

\noindent where $\boldsymbol{\varphi}^i$ and $\boldsymbol{x}^i$ denote the deformation map and the deformed configuration after the $i^{th}$ Newton step. According to (10), the deformation is repeatedly reset to the initial undeformed state to erase the complete deformation history within the fictitious domain $\Omega_{\textit{fict}}$. This does not affect the physical consistency and accuracy of the solution in the physical domain $\Omega_{\textit{phys}}$, provided that the influence of $\Omega_{\textit{fict}}$ is mitigated by a sufficiently strong penalization. Furthermore, the assumption of (10) supersedes the calculation of the deformation gradient \cite{Schillinger:12.2}, so that any stability issues resulting from the numerical computation of the deformation gradient are automatically avoided. The corresponding stress solution is plotted in Fig.~8c, where the oscillatory behavior of Fig.~8b is considerably reduced by several orders of magnitude. Moreover, the deformation resetting can be efficiently implemented by exploiting the coincidence of linear and geometrically nonlinear elasticity at the deformation and stress-free reference configuration \cite{Schillinger:12.2}.

To test convergence in an energy measure, the uni-axial rod of Fig.~5 is considered with sine-load $f_{\text{sin}}$ and $\Delta u$=1.0. The convergence under $p$-refinement is plotted in Fig.~9 for the $p$- and B-spline versions. For the standard FCM formulation, it illustrates the convergence decay to a low algebraic rate as a consequence of the insufficient penalization in conjunction with oscillatory stresses. The modified geometrically nonlinear formulation allows for a decrease of the penalty parameter to $\alpha$=$10^{-15}$, which restores the ability of the FCM schemes to achieve exponential convergence.

\section{Extension to Transport Processes in Porous Media}

The concept of the finite cell method has also been applied to other types of partial differential equations, e.g.~to transport problems within the MAC project B5 \textit{Transport and Reaction Processes in Porous Media}. From a geometric point of view, simulation of transport processes in porous media may be similarly challenging as the structural problems described in the previous sections. The domain of computation is often very complex, and the generation of a mesh that conforms to the porous structure may be highly involved.

We start from the weak formulation of the transport equation that reads

\begin{equation}
\int_{\Omega}[\vec{q} c \cdot\nabla w + (\theta\nu) \nabla c \cdot\nabla w]d\Omega = \int_{\Omega}w f d\Omega
\end{equation}

\noindent where $c$ is the concentration, $w$ is a test function and $f$ is a source term. The quantities $\vec{q}$ and $(\theta\nu)$ denote the Darcy's velocity and the effective diffusion coefficient, respectively \cite{Helmig:97.1}. In the sense of Fig.~1, the porous flow domain $\Omega_{\textit{phys}}$ is embedded in a larger domain, which is meshed by a simple Cartesian grid. The bilinear form associated to the weak form of (11) is extended to the fictitious domain, and the quantities $\vec{q}$ and $(\theta\nu)$ are multiplied by a penalty factor $\alpha$ in the sense of (3):
\begin{eqnarray}
\left\{ \begin{array}{ll}
\vec{q}_e = \alpha\cdot\vec{q} \\
(\theta\nu)_e = \alpha\cdot\theta\nu
\end{array}\right.
\label{eqn:qe}
\end{eqnarray}

\noindent Quadtree-based numerical integration is performed for cells cut by the boundary of the flow domain (see Section~2.3) and a Bubnov-Galerkin Ansatz is made with high-order basis functions\footnote{It is worthwhile to note that high-order basis functions are significantly more stable than low-order functions for flow problems moderately dominated by convection \cite{Cai:13.1}.}.

\begin{figure}[t!]
\centering
\subfloat[Concentration profile for a P\'eclet number of Pe=1 \cite{Donea:03.1}.] {\hspace{3cm}\includegraphics[width=0.62\textwidth]{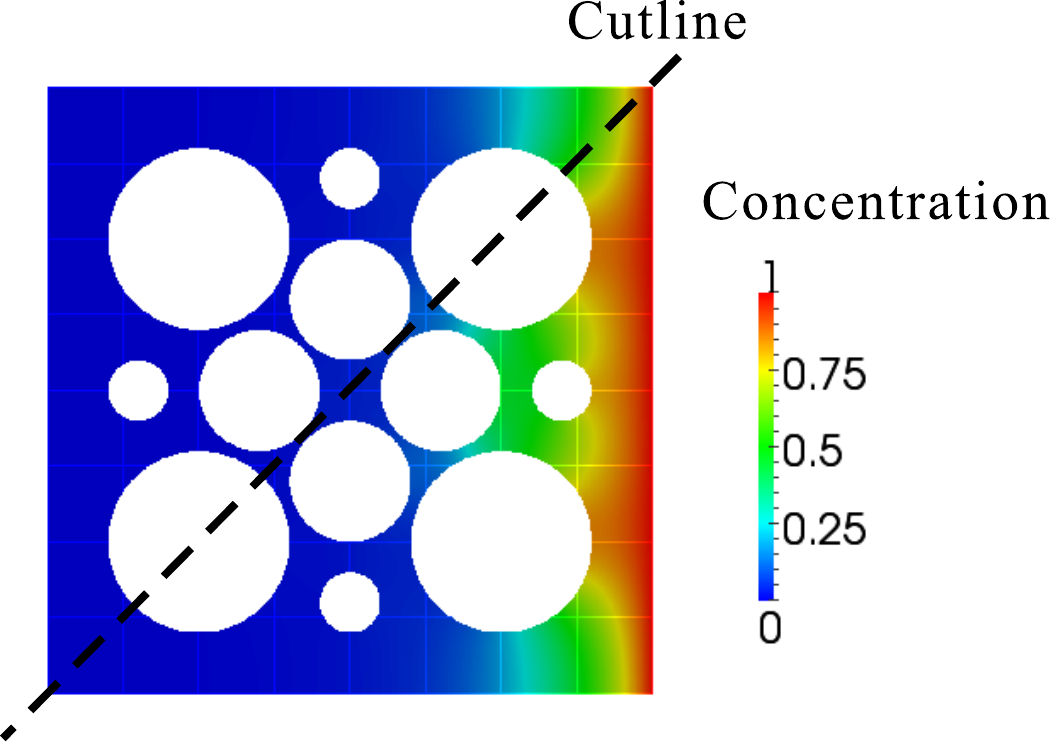}\hspace{3cm}} %\hspace{0.69cm}

\subfloat[Finite cell solution along the diagonal cut line as compared to a body-fitted $p$-FEM reference solution.]{\hspace{2cm}\includegraphics[width=0.62\textwidth]{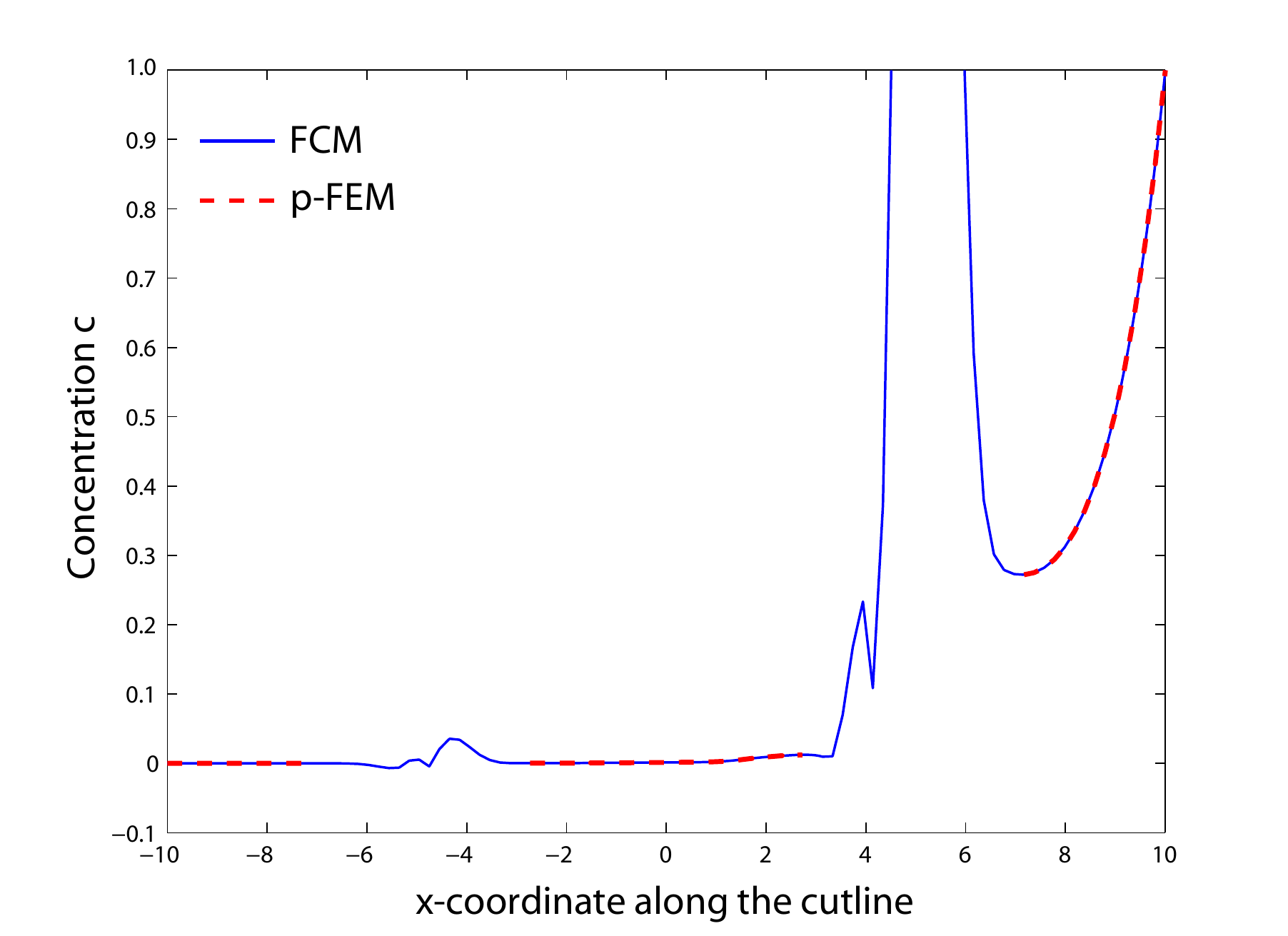}\hspace{2cm}}
\caption{FCM analysis of porous media: The test case consists of a square domain with impermeable spheres. Dirichlet constraints of c=0 and c=1 are imposed on the left and right boundaries, respectively, and no flow conditions on the upper and lower boundaries.}
\end{figure}

Figure~10 shows results of a test case for the simulation of transport through porous media, where for simplicity the velocity $\vec{q}$ of the transporting fluid was assumed to be constant along the x-axis throughout the domain. Although this can only be assumed for a limited range of physical problems, it does not impose a restriction with respect to the validation of the approach, as long as the transport velocity $\vec{q}$ is taken as an upper bound of the expected true flow velocity. The square domain is discretized by $8\times8$ finite cells of polynomial degree $p$=8. Figure~10b illustrates the profile of concentration $c$ along the diagonal cut line. A reference solution obtained by a refined computation of a boundary fitted mesh is compared to the FCM solution for a P\'eclet number of Pe=1 \cite{Donea:03.1,Helmig:97.1}. Similar to structural problems, the FCM solution shows very good quality, even in areas, where only ``narrow" flow bridges between obstacles are present.

\begin{figure}[t!]
\centering
\subfloat[B\'{e}zier elements of a T-spline surface (Output from CAD package Rhino with  T-spline plug-in).]{\hspace{3cm}\includegraphics[width=0.38\textwidth]{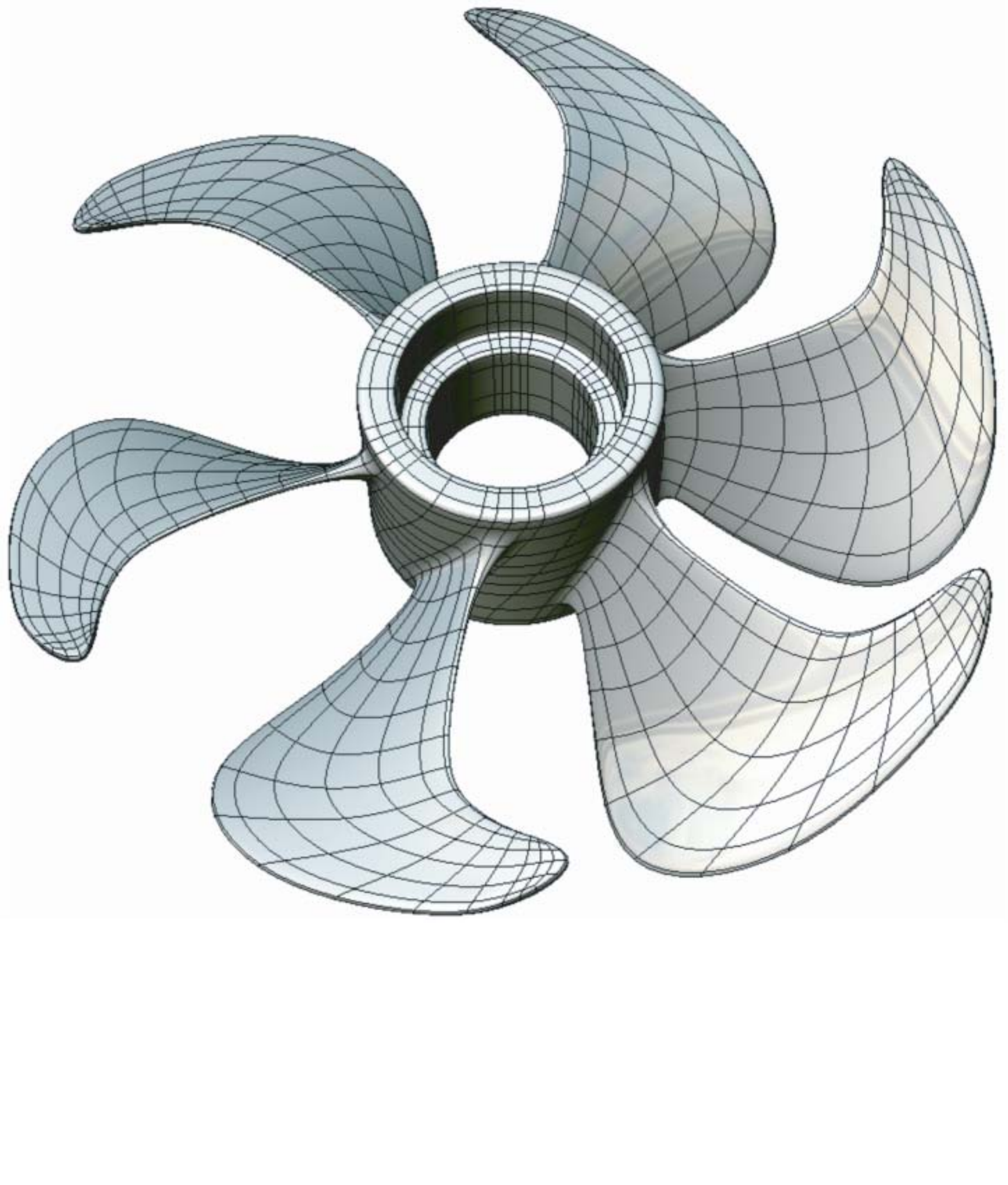}\hspace{3cm}} %\hspace{0.69cm}

\subfloat[The complete structure is immersed in a bounding box of 16$\times$16$\times$4 axis-aligned cubic B-spline elements.]{\hspace{2cm}\includegraphics[width=0.55\textwidth]{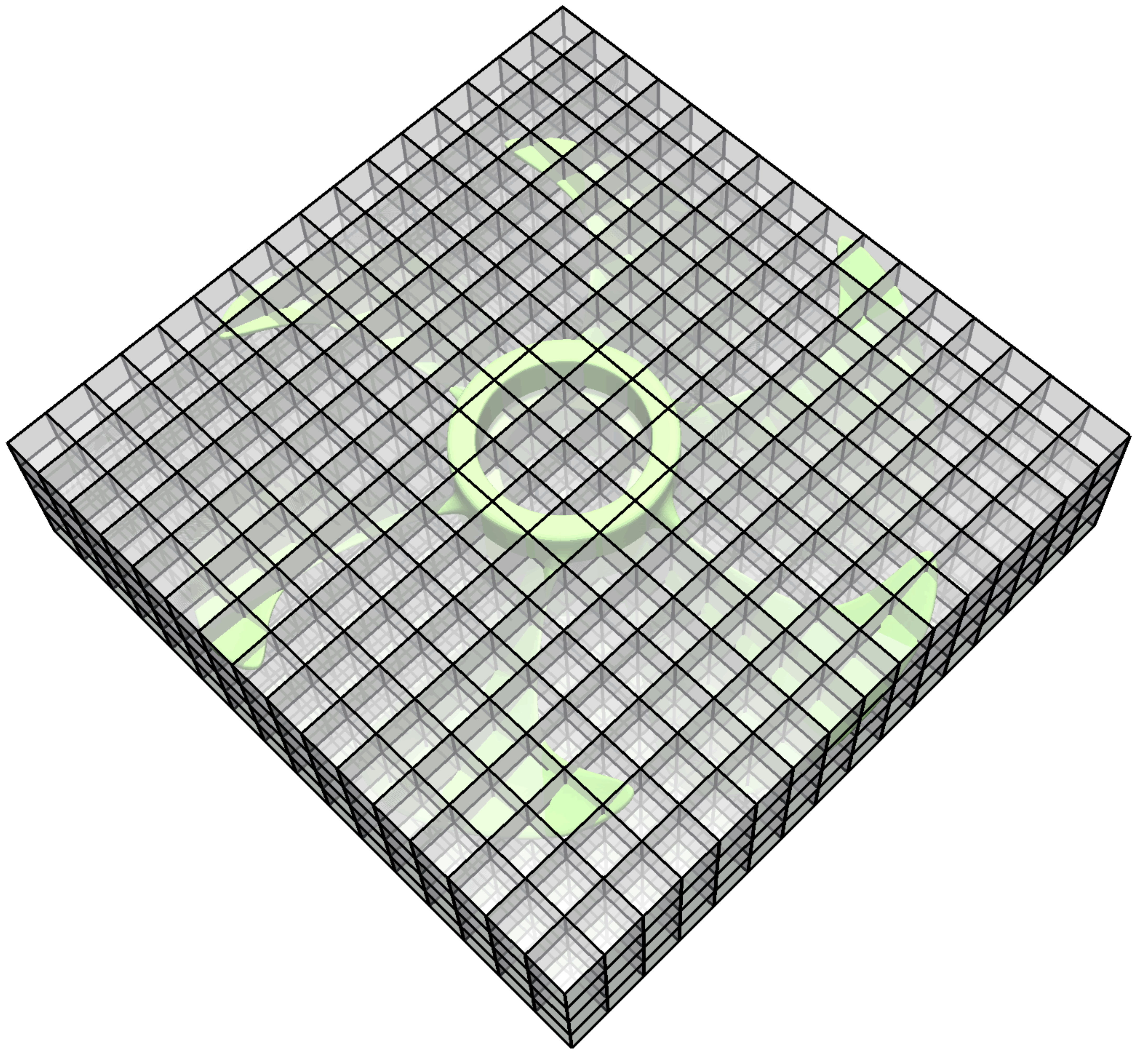}\hspace{2cm}}
\caption{Ship propeller example: CAD based geometry description and finite cell discretization.}
\end{figure}

\begin{figure}[t!]
\centering
\subfloat[Deletion of elements without support in the propeller domain creates a reduced set of elements, which homogeneously resolve the structure irrespective of the local thickness.]{\hspace{3.6cm}\includegraphics[width=0.476\textwidth]{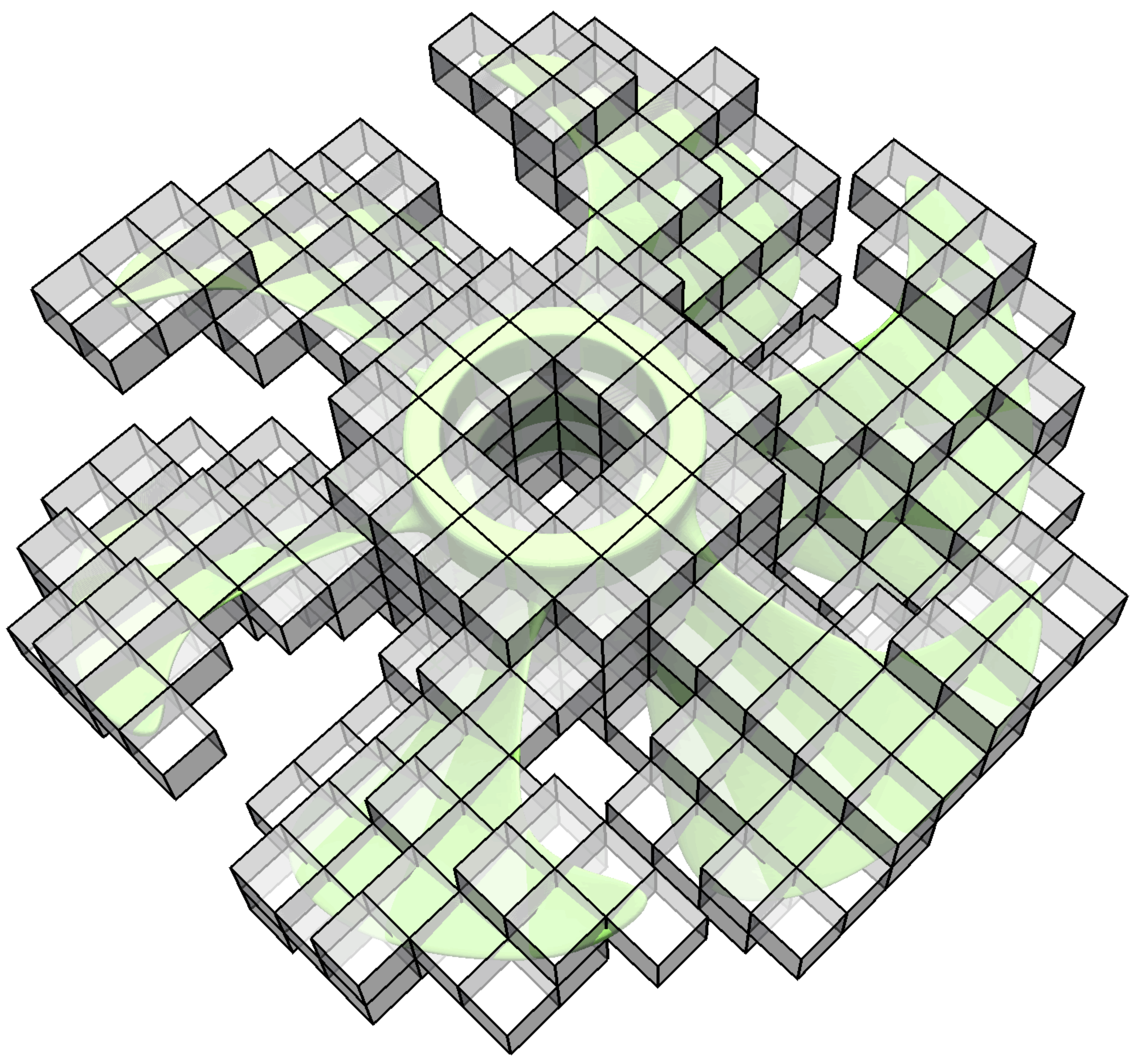}\hspace{3.6cm}} 

\subfloat[Hierarchical refinement of the propeller blades achieves a homogeneous through-the-thickness resolution.]{\hspace{3.6cm}\includegraphics[width=0.476\textwidth]{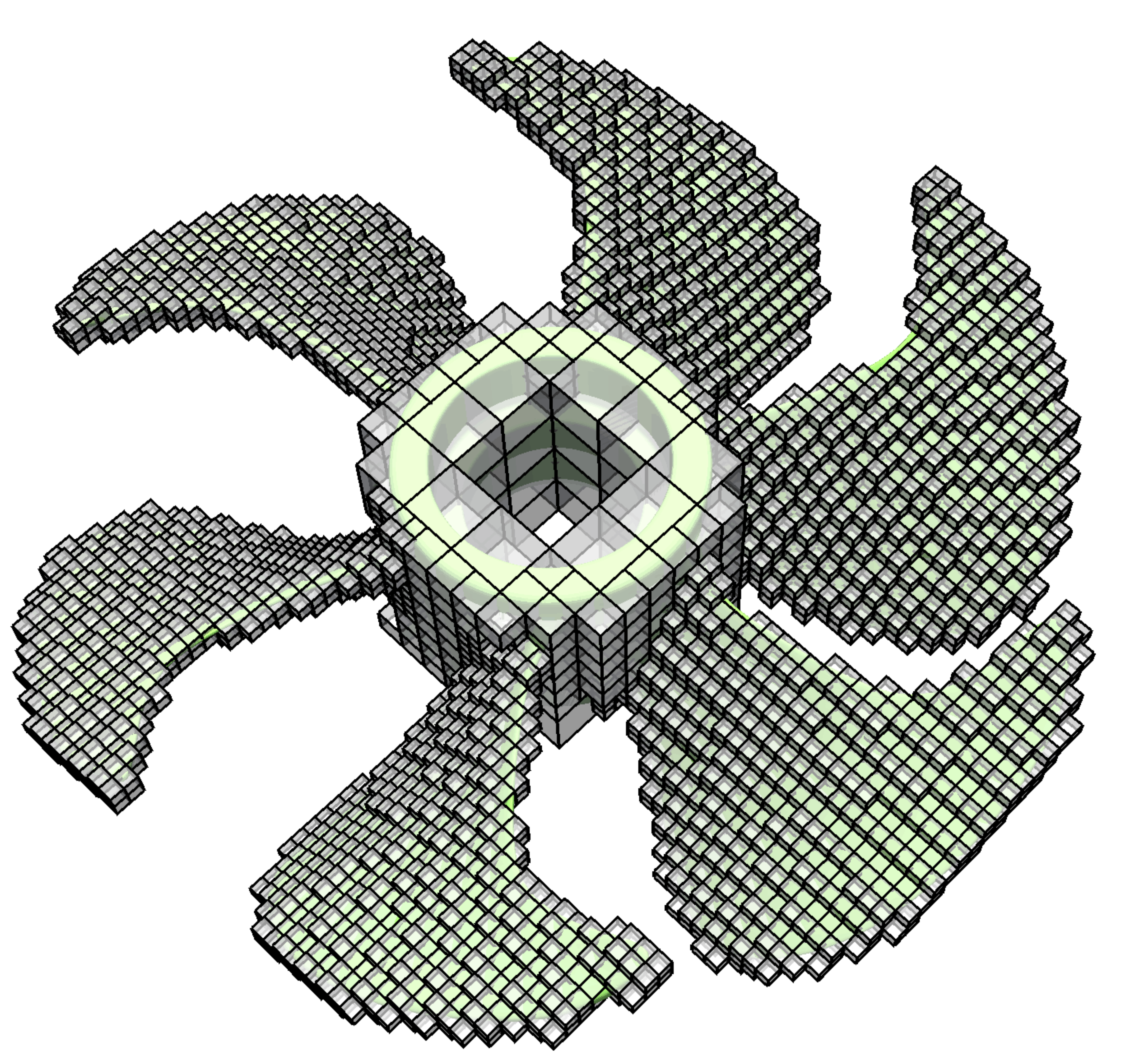}\hspace{3.6cm}}
\caption{Ship propeller example: The role of hierarchical refinement.}
\end{figure}

\section{Application Oriented Examples: Structural Analysis of CAD and Image-based Geometric Models}

In the following, we illustrate the benefits of the finite cell method in terms of simple mesh generation for very complex geometries by two application oriented examples, which are described by a CAD (computer aided design) based T-spline surface and a CT (computed tomography) based voxel model, respectively. A more detailed description and further computations for the example structures can be found in \cite{Schillinger:12.2,Schillinger:12.1,Schillinger:12.3}.

\subsection{Modal Analysis of a Ship Propeller}

The geometry of the propeller is given by a smooth, watertight T-spline surface (i.e., there are no gaps or overlaps). It is exported from the CAD package Rhino \cite{Rhino:12.1} in conjunction with the T-spline plug-in in the form of B\'ezier elements as shown in Fig.~11a. Its maximum diameter and height is 0.695 m and 0.334 m, respectively, and it is made out of steel with Young's modulus 2.1$\cdot 10^{11}$ N/m$^2$, Poisson's ratio 0.28 and density 7,850 kg/m$^3$. The structure can neither be characterized as a typical shell nor as a true solid. Configurations like this usually require specialized and time consuming meshing procedures to produce good quality discretizations.

\begin{figure}[t!]
\centering
\includegraphics[width=0.476\textwidth]{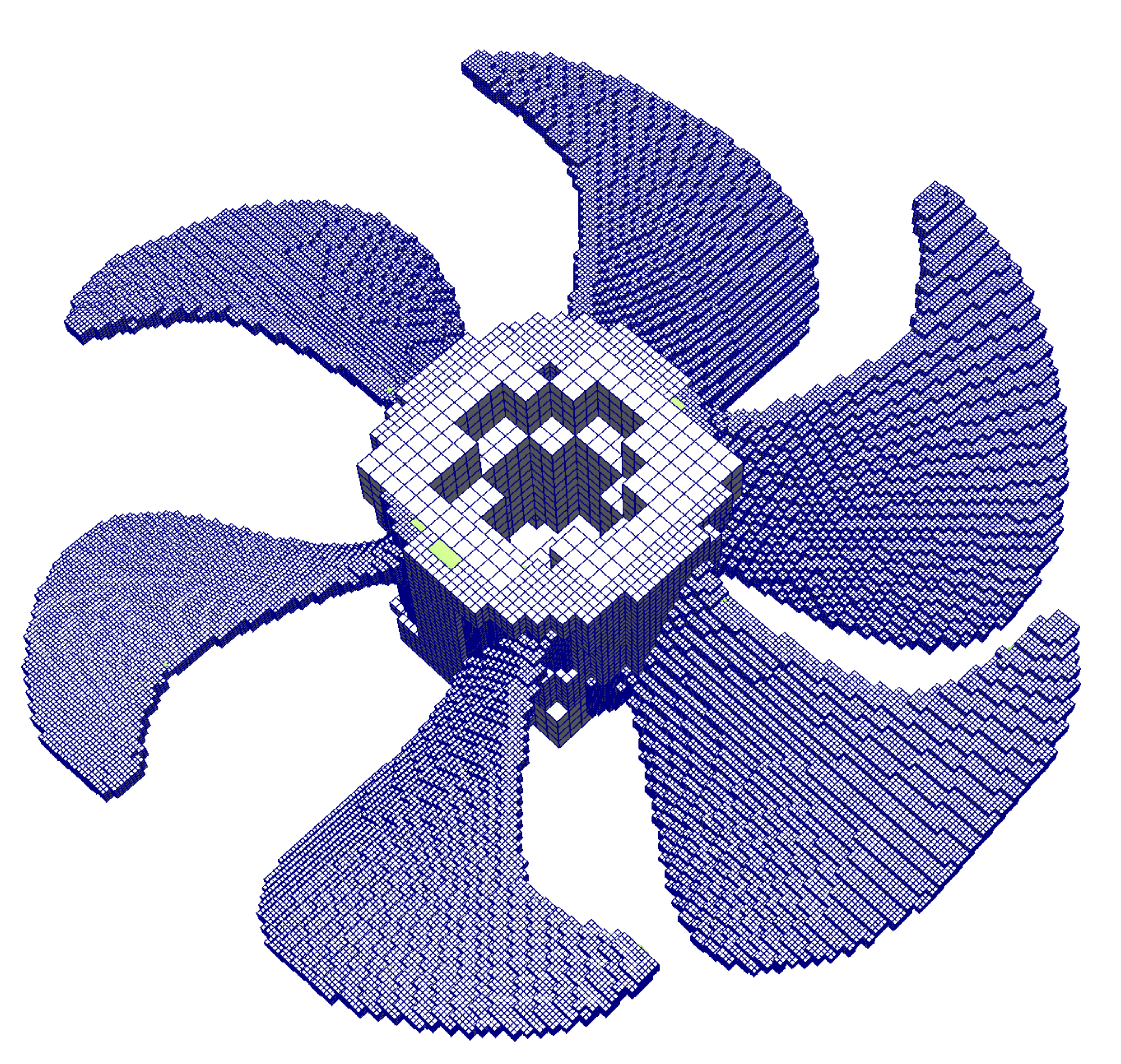}
\caption{Sub-cell partitioning of elements cut by the geometric boundary. The adaptive decomposition scheme shown in Fig.~4 is carried out up to level $k$=2.}

\includegraphics[width=0.49\textwidth]{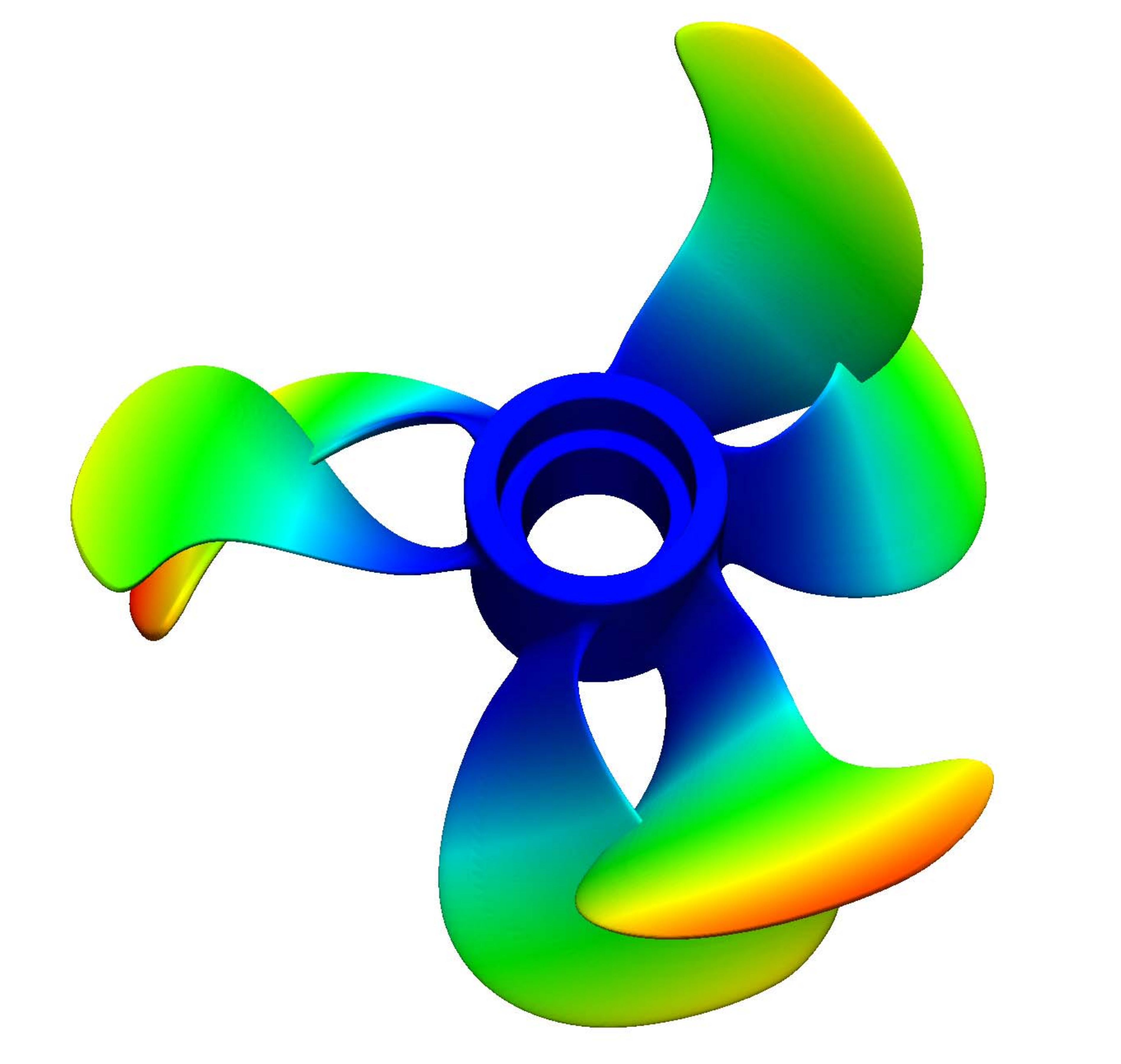}
\caption{The first, second and third mode shapes exhibit a rotational symmetry around the center, corresponding to the three pairs of opposing propeller blades. We display mode 2.}
\end{figure}

The finite cell method circumvents the meshing challenge completely, which we demonstrate in the following with the B-spline version of the FCM. First, the complete structure is embedded in a regular grid of axis-aligned B-splines of polynomial degree $p$=3 (see Fig.~11b). Second, all knot span cells without support in the propeller domain are eliminated from the discretization (see Fig.~12a). The decision whether an element is to be kept or not is based on a simple point location query, which checks if at least one integration point is located in $\Omega_{\textit{phys}}$. It can be efficiently implemented for example by search algorithms based on special space-partitioning data structures such as $k$-d trees \cite{Samet:06.1,Bindick:09.1}. An axis-aligned discretization with elements of the same size does not account for the inhomogeneous thickness of the different regions of the structure. In a third step, we therefore apply two levels of hierarchical refinement to the propeller blades, while we leave the discretization of the central hub as it is, to achieve a homogeneous resolution of the two different thicknesses \cite{Schillinger:12.3}. In a fourth step, we equip each element cut by the geometric boundary by additional sub-cells, which are organized in an octree of depth two (see Fig.~13). Each sub-cell contains 4$\times$4$\times$4 Gauss points, leading to an aggregation of integration points in cut elements to accurately take into account the geometric boundary during numerical integration. The contribution to stiffness and mass matrices that result from integration points located outside the propeller domain $\Omega_{\textit{phys}}$ are penalized by factor $\alpha$=10$^{-3}$. The hierarchically refined mesh of Fig.~12b is analysis suitable and is used in combination with the sub-cells of Fig.~13 to conduct a modal analysis of the structure, where the mass matrix is lumped according to the row sum method \cite{Hughes:00.1}. Figure~14 illustrates the first mode shapes.

\begin{figure}[t!]
\centering
\subfloat[Voxelized sample cube.] 
{\includegraphics[width=0.5\textwidth]{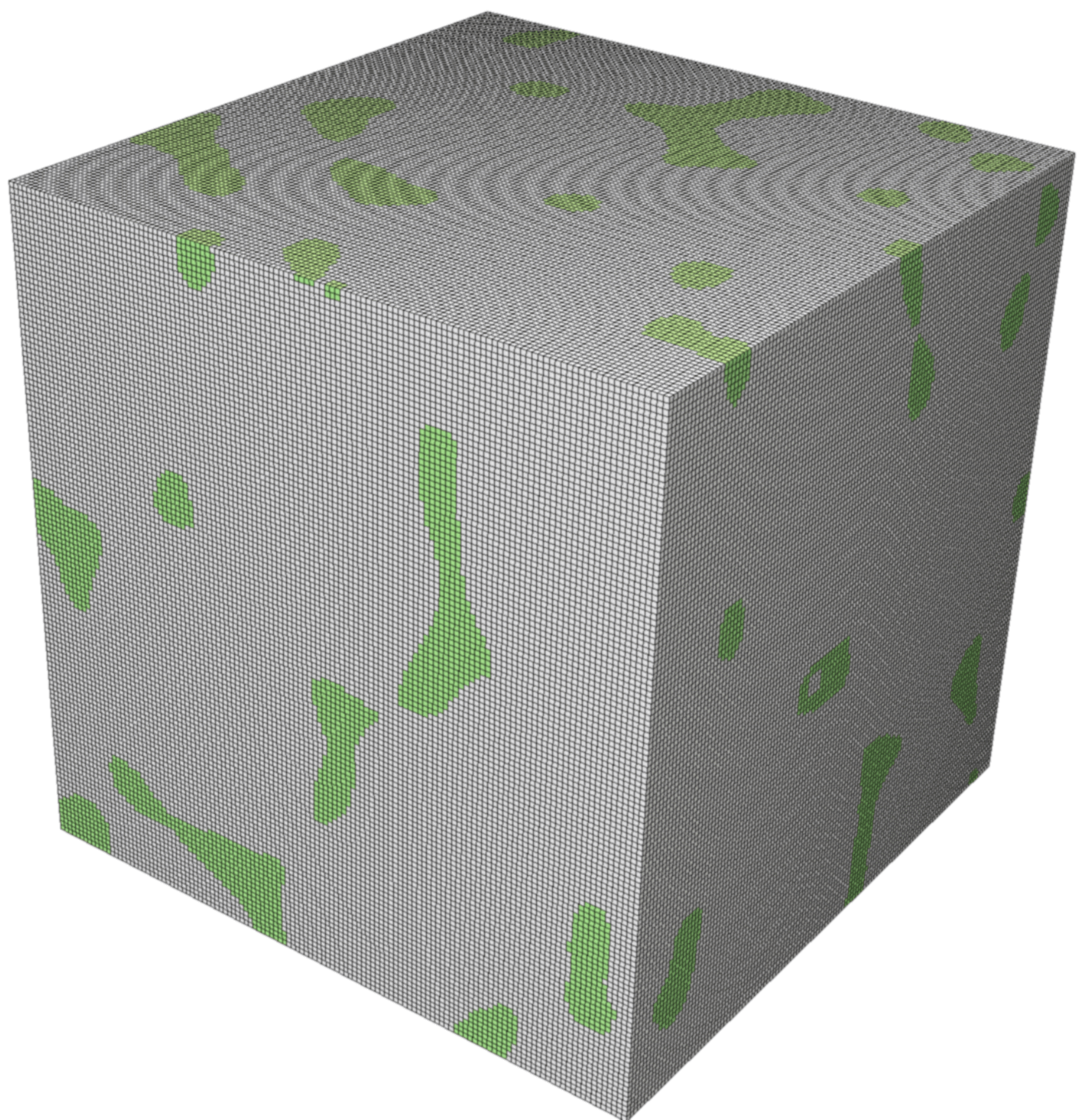}}

\subfloat[Voxels with $b_{\textit{vox}}=1$.] 
{\includegraphics[width=0.5\textwidth]{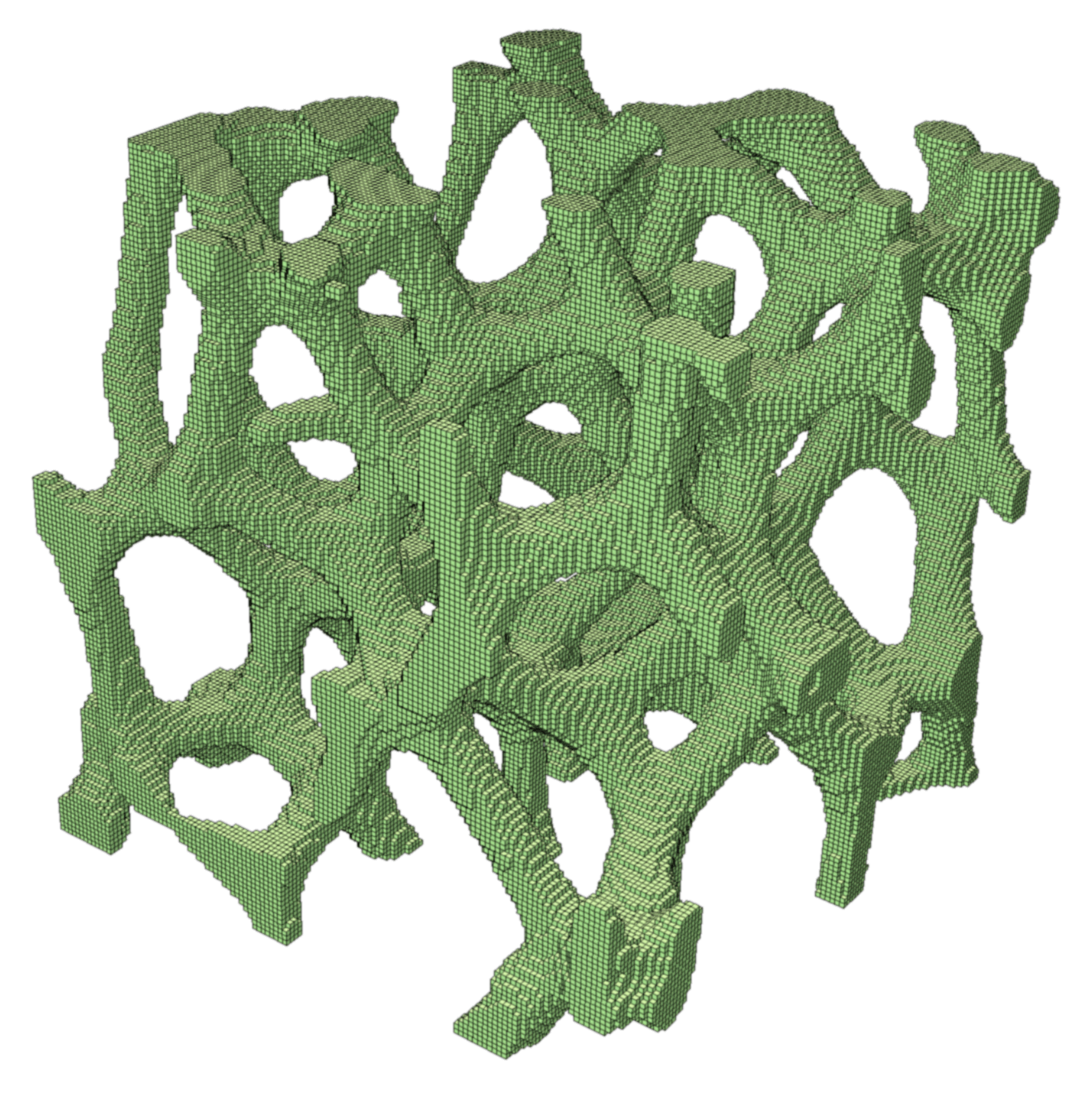}}
\caption{CT-based voxel model of an aluminium foam sample. For better visibility, the original resolution of 1024$^3$ voxels is reduced to 128$^3$.}
\end{figure}

\begin{figure}[t!]
\centering
\includegraphics[width=0.6\textwidth]{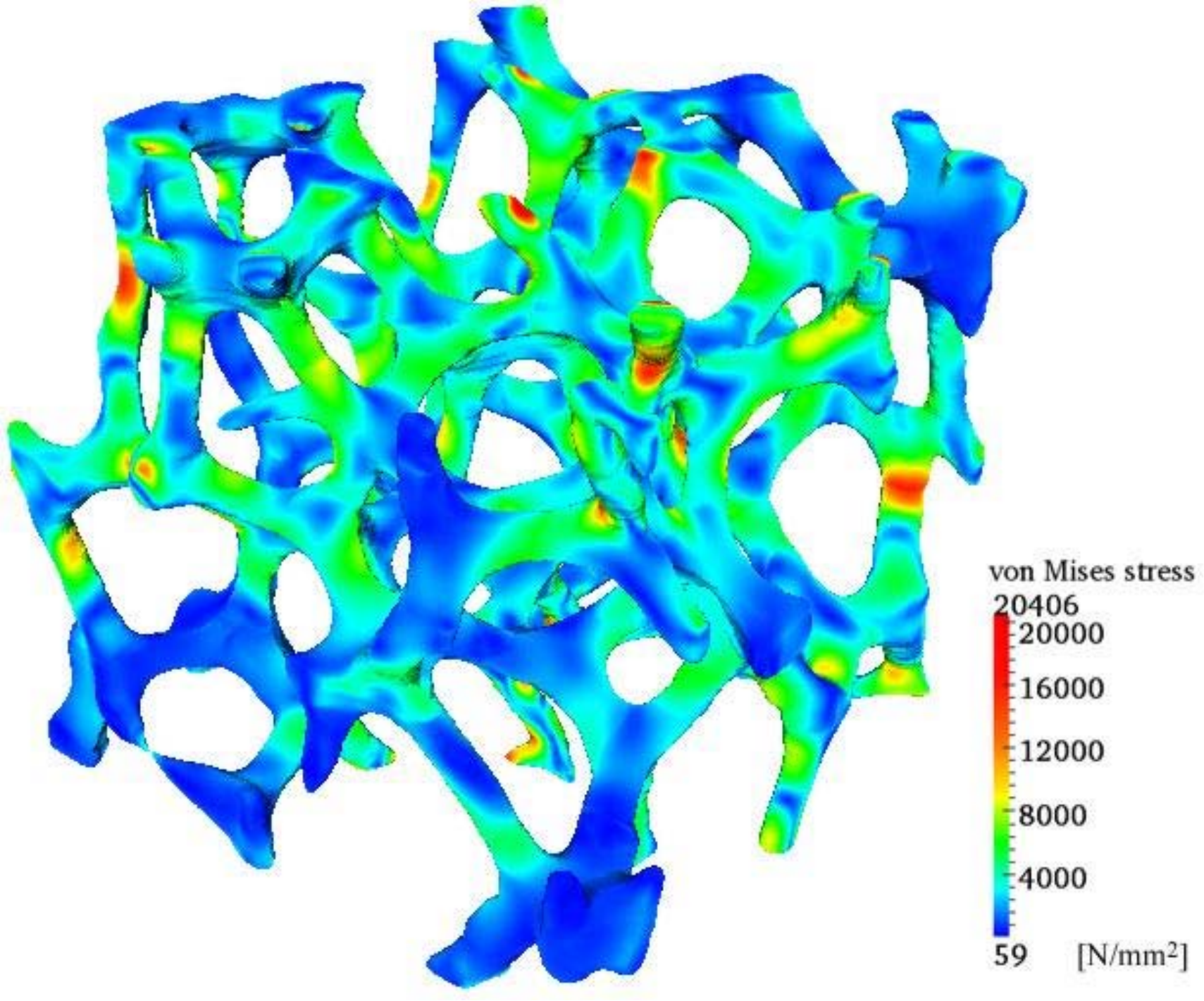}
\caption{Von Mises stress of the metal foam sample plotted on the deformed configuration. The results are obtained from a discretizations with $5\times5\times5$ finite cells of polynomial degree $p$=7.}
\end{figure}

\subsection{Large Deformation Analysis of an Open-cell Aluminium Foam}

Metal foams provide high stiffness at reduced weights, and are therefore frequently used for lightweight structures in automotive and aerospace applications \cite{Banhart:01.1}. The $p$-version of the FCM is applied to simulate a compression test for an aluminium foam sample of size 20$\times$20$\times$20 mm, discretized by a structured grid of $5\times5\times5$ high-order finite cells. Its internal geometry is provided by voxels with a resolution of $1024$ in each Cartesian direction, each of which encodes $\alpha$. Figures~15a and 15b show the complete voxel model of the sample cube and the physical voxels of material index 1 associated with aluminium, respectively, in a coarsened resolution of $128^3$. The foam sample is assumed as part of a larger specimen, which is uniformly compressed along the vertical axis. Corresponding boundary conditions are specified as follows \cite{Sehlhorst:09.1}: Displacements normal to the top surface are gradually increased to 1.6mm (8\% compressive deformation), modelling the influence of a testing machine, whereas the displacements normal to all other surfaces are fixed due to the bottom support and the influence of the surrounding material of the specimen. The aluminium foam is characterized by Young's modulus \textit{E}=70.000 N/$mm^2$, penalized by $\alpha$=$10^{-12}$ at all integration points in $\Omega_{\textit{fict}}$, and Poisson's ratio $\nu$=$0.35$.

The finite cell method is able to directly operate on the voxel model, which provides a basis for a simple point location query (see \cite{Schillinger:12.2}). In particular, we avoid the costly transformation of the voxel model into a surface model by image-based software, which is required as a basis for mesh generation by standard body-fitted simulation methods. For the present $p$-version mesh of polynomial degree $p$=7 with 3 levels of sub-cells (21,492 dofs; 24,947 sub-cells; approx. 12.75 million Gauss points), analysis of the foam could be accomplished by our in-house FCM code in about 4 hours\footnote[1]{Using 8 threads on 2 interconnected Intel(R) Xeon(R) W5590 @ 3.33GHz}. Since the major cost of FCM results from the large number of sub-cells with full Gauss integration, a major performance gain is achieved by the shared memory parallelization of the loop that computes local stiffness matrices for cells and sub-cells with subsequent assembly into the global system matrix. A \textit{parallel for} construct creates a team of $n$ threads to execute the main loop over sub-cells in parallel, where $n$ is the number of threads available. With $n$=8, we achieved a strong speed up of the loop of around 5. In addition, we would like to mention here that a very fast variant of the FCM has been developed (see \cite{Knezevic:11.1,Knezevic:11.2,Yang:12.1,Yang:12.2}), relying on pre-integration of sub-element matrices. This implementation is up to three orders of magnitude faster than the FCM based on standard quadrature described herein, and thus allows even for real-time simulations of complex 3D structures.

The resulting von Mises stresses shown in Fig.~16 exhibit accurate localization of stress concentrations at the convex sides of the foam members, which agrees well with engineering experience. Figure~17 plots the equivalent force obtained from integration of the normal stress over the top surface vs. the prescribed displacement of the top surface for different polynomial degrees $p$. It can be observed that the increase of $p$ improves the reproduction of the geometrically nonlinear behavior of the foam.
  
\begin{figure}[t!]
\centering
\includegraphics[width=0.6\textwidth]{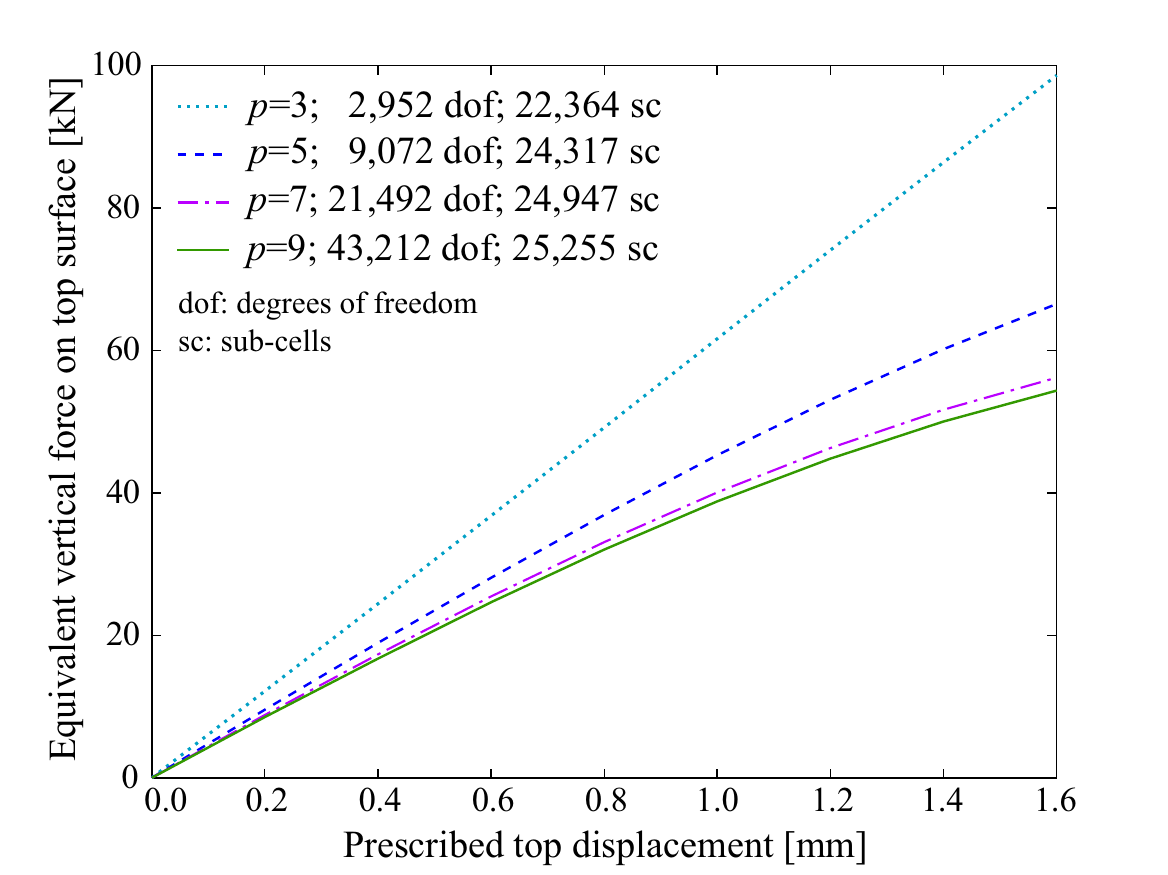}
\caption{Convergence of the force-displacement behavior under $p$-refinement for the foam sample.}
\end{figure}

\section{Summary and Outlook}

The present contribution provided a review of recent developments of the finite cell method (FCM) for the analysis of complex structures. We briefly summarized the basic components of the FCM technology, i.e.~the fictitious domain concept, high-order basis functions, adaptive integration and weak boundary conditions, and outlined its basic numerical properties for linear elasticity, i.e.~the smooth extension of solution fields beyond the physical domain as well as exponential rates of convergence in energy norm. We then summarized the concept of deformation resetting, which enables the extension of the FCM to nonlinear elasticity. Finally, we illustrated with two application examples, i.e.~a ship propeller and a metal foam sample, that the benefits of the finite cell method in terms of almost no mesh generation for complex structures can be achieved for both CAD based and explicit image-based geometric models.

Based on these results, we believe that the finite cell method has great potential for the accurate analysis of very complex structures, and a plethora of very promising aspects are still open, such as the analysis of topology changes and moving boundaries, for which embedded domain methods such as the FCM offer significant advantages over ALE-type approaches, or the introduction of FCM suitable coupling schemes for multiphysics problems, which stand at the forefront of today's challenges in computational science.\\[0.36cm]

\noindent $\bf{Acknowledgments.}$ D.~Schillinger, Q.~Cai and R.-P.~Mundani gratefully acknowledge support from the Munich Centre of Advanced Computing (MAC) and the International Graduate School of Science and Engineering (IGSSE) at the Technische Universit\"{a}t M\"{u}nchen. D.~Schillinger gratefully acknowledges support from the German National Science Foundation (Deutsche Forschungsgemeinschaft DFG) under grant number SCHI 1249/1-1. The authors thank T.J.R.~Hughes, M.~Ruess and M.A.~Scott for their help with the analysis of the ship propeller.\\

\bibliographystyle{unsrt}
\bibliography{paper}

\end{document}